\documentclass[a4paper,11pt]{article}

\usepackage{microtype}
\usepackage[]{algorithm2e}
\usepackage{chngpage}

\usepackage[T1]{fontenc}
\usepackage[utf8]{inputenc}

\usepackage{amsfonts}
\usepackage{amsmath}
\numberwithin{equation}{section}

\usepackage{amssymb}
\usepackage{mathtools}
\usepackage[binary-units=true]{siunitx}

\usepackage{bm}
\usepackage{breqn}

\usepackage{amsthm}
\newtheorem{theorem}{Theorem} 
\newtheorem{lemma}[theorem]{Lemma}

\usepackage{caption}
\usepackage{subcaption}
\usepackage{diagbox}

\usepackage{txfonts}
\usepackage{sansmath}

\usepackage{geometry}
\newgeometry{vmargin={25mm}, hmargin={25mm,25mm}}

\usepackage{multicol}
\usepackage{booktabs}
\usepackage{graphicx}
\usepackage{float}

\makeatletter
\newcommand\tabfill[1]{%
\dimen@\linewidth%
\advance\dimen@\@totalleftmargin%
\advance\dimen@-\dimen\@curtab%
\parbox[t]\dimen@{#1\ifhmode\strut\fi}%
}

\usepackage{dcolumn}
\newcolumntype{d}[1]{D{.}{.}{#1}}

\usepackage{tikz}
\usetikzlibrary{shapes,arrows}
\usetikzlibrary{matrix}
\usepackage{pgfplots}
\usepackage{pgfplotstable}
								
\usepackage{color}
\definecolor{thered}{rgb}{0.65,0.04,0.07}
\definecolor{thegreen}{rgb}{0.06,0.44,0.08}
\definecolor{theblue}{rgb}{0.02,0.2,0.68}
\definecolor{clr1}{HTML}{bdd7e7}
\definecolor{clr2}{HTML}{6baed6}
\definecolor{clr3}{HTML}{3182bd}
\definecolor{clr4}{HTML}{08519c}

\newtheorem{remark}{Remark}
\newtheorem{corollary}{Corollary}[theorem]

\renewcommand{\@algocf@capt@plain}{above}

\usepackage[unicode=true,
            bookmarks=true,
            bookmarksnumbered=false,
            bookmarksopen=false,
            breaklinks=true,
            pdfborder={0 0 1},
            backref=false,
            colorlinks=true,
            linkcolor=black,
            urlcolor=theblue,
            citecolor=theblue,
            hyperfootnotes=false]
           {hyperref}
\hypersetup{pdftitle={},
            pdfauthor={B. C. Vermeire}}

\begin{document}

\title{Quasi-Orthogonal Runge-Kutta Projection Methods}

\author{Mohammad R. Najafian and Brian C. Vermeire\\\\
  \textit{Department of Mechanical, Industrial, and Aerospace Engineering}\\
	\textit{Concordia University}\\
	\textit{Montreal, QC, Canada}}
\maketitle

\begin{abstract}
A wide range of physical phenomena exhibit auxiliary admissibility criteria, such as conservation
of entropy or various energies, which arise implicitly under exact solution of their governing PDEs. However, standard temporal schemes, such as classical Runge-Kutta (RK) methods, do not enforce these constraints, leading to a loss of accuracy and stability. Projection is an efficient way to address this shortcoming by correcting the RK solution at the end of each time step. Here we introduce a novel projection method for explicit RK schemes, called a \textit{quasi-orthogonal} projection method. This method can be employed for systems containing a single (not necessarily convex) invariant functional, for dissipative systems, and for the systems containing multiple invariants. It works by projecting the orthogonal search direction(s) into the solution space spanned by the RK stage derivatives. With this approach linear invariants of the problem are preserved, the time step size remains fixed, additional computational cost is minimal, and these optimal search direction(s) preserve the order of accuracy of the base RK method. This presents significant advantages over existing projection methods. Numerical results demonstrate that these properties are observed in practice for a range of applications.
\end{abstract}

\section{Introduction}

Many physical systems yield measurable quantities that either remain constant or evolve monotonically over time. Examples of such behavior include explicitly conserved quantities, such as mass, momentum, and energy, or implicitly conserved quantities like entropy as in the compressible Euler equations with smooth solutions. When solving these Partial Differential Equations (PDEs) numerically, an important indicator of the quality of the numerical solution is preserving physical progression of these parameters \cite{hairerGeometricNumericalIntegration2006}. Often, these constraints are enforced explicitly through chosen conservation laws, such as conservation of mass, momentum, and total energy. However, in some cases such as entropy, they are enforced indirectly under the exact evolution of conservation laws. Failure to maintain these with approximate discrete solutions can result in non-physical solutions, stability issues, or increased numerical error for long-time integration
\cite{gearInvariantsNumericalMethods1992, arakawaComputationalDesignLongTerm1997, calvoErrorGrowthNumerical2011, ketchesonRelaxationRungeKutta2019, ranochaRelaxationRungeKutta2020, ranochaRateErrorGrowth2021}. Therefore, numerical schemes that can provide structure-preserving solutions, maintaining these implicitly conserved quantities, are of great importance.
  
Solving time-dependent PDEs commonly involves two main steps: first, the domain undergoes spatial discretization, transforming the system into semi-discrete time-dependent Ordinary Differential Equations (ODEs). Subsequently, these ODEs are solved using a temporal integration scheme to obtain a numerical solution that is fully-discrete in space and time. Each of these steps can potentially contaminate conserved quantities. For spatial discretization, there is a rich literature on stability-preserving techniques, for example \cite{delreyfernandezReviewSummationbypartsOperators2014,  chanDiscretelyEntropyConservative2018, chenReviewEntropyStable2020, kuyaHighorderAccurateKineticenergy2021} for Euler and Navier-Stokes equations. Nevertheless, the resulting semi-discrete ODEs should be coupled with a time integration scheme to proceed in time, and the important question remains: if the conservation of these quantities is maintained after numerical time integration. To answer this, stand-alone integration schemes can be studied in terms of conserving the nonlinear stability properties of the ODE systems. 

To proceed, consider the following time-dependent ODE system, which may represent a PDE problem after spatial discretization 
\begin{subequations}\label{IVP}
  \begin{align}
    & \frac{d}{dt}\underline{q} \left(t \right) = R \left( t, \underline{q} \left( t \right) \right) ,\\
    & \underline{q} \left(0 \right) = \underline{q}_0 ,
  \end{align}
\end{subequations} 
where $\underline{q} : \mathbb{R} \rightarrow \mathbb{R}^m$  and $R : \mathbb{R} \times \mathbb{R}^m \rightarrow \mathbb{R}^m  $. Suppose that for this system there is a  smooth nonlinear function $G \left(\underline{q} \right) : \mathbb{R}^m \rightarrow \mathbb{R}$, called an invariant, whose time derivative is zero for the time range of interest
\begin{equation} \label{dG_dt}
\frac{d}{dt} G \left( \underline{q} \left(t \right) \right)=  \nabla G \left(\underline{q} \right)^T R \left(t,\underline{q} \right)= 0,  \qquad \textrm{for all} \; t \in [0, T].
\end{equation}
Preserving this property by an integration scheme implies that at each time step the invariant magnitude remains constant up to machine precision. However, many widely used integration schemes, including Runge-Kutta (RK) schemes, cannot guarantee conservation of general nonlinear invariants \cite{hairerGeometricNumericalIntegration2006}.

With an standard $s$-stage RK scheme and having the solution at the current time step, $\underline{q}^{n}$, the approximate next step solution, $\underline{q}^{n+1}$, is obtained as

\begin{subequations}
\begin{align}
& \underline{q} \left(t_n + \Delta t \right) \simeq \underline{q}^{n+1} = \underline{q}^n + \Delta t \sum_{j=1}^s b_j\underline{R}_{j}, \\
& \underline{R}_{i} = R \left(t_n + c_i \Delta t, \underline{q}_i \right),\quad \underline{q}_i= \underline{q}^n + \Delta t \sum_{j=1}^s a_{ij}\underline{R}_{j}, \quad c_i= \sum_{j=1}^s a_{ij},
\end{align}
\end{subequations}
where $a_{ij}$ and $b_j$ are coefficients of the selected RK method and $\underline{R}_{j} $ represents stage derivatives. These derivatives can be calculated explicitly if the RK scheme is explicit, i.e. $a_{ij}=0$ for $j \geq i$. 

All standard RK schemes automatically preserve linear invariants, e.g. total mass. Assume that $L\left(\underline{q}\right)= \underline{d}^T \underline{q}$ with a constant vector $\underline{d}$ is a linear invariant of the problem. It means that $d^TR\left(t, \underline{q}\right)$ will be zero for all $t\in [0, T]$ and $\underline{q} \in \mathbb{R}^m$. So, with a RK method every stage derivative $\underline{R}_j$ is orthogonal to $\underline{d}$
\begin{equation}
\underline{d}^T\underline{R}_j=0.
\end{equation}
In other words, moving the solution along stage derivatives will not change linear invariants of the system. Consequently, with standard RK integration schemes linear invariants of the problem will be preserved after each time step \cite{hairerGeometricNumericalIntegration2006}
\begin{equation}
 \underline{d}^T \underline{q}^{n+1}= \underline{d}^T\underline{q}^n + \underline{d}^T \Delta t \sum_{j=1}^s b_j\underline{R}_j= \underline{d}^T \underline{q}^n. 
\end{equation}

However, no explicit RK scheme can guarantee conservation of general quadratic invariants, and neither explicit nor implicit RK schemes can conserve general nonlinear invariants of order three or higher \cite{hairerGeometricNumericalIntegration2006}. In fact, RK schemes only preserve nonlinear invariants up to truncation error, not up to machine precision 
\begin{equation}
G\left(\underline{q}^{n+1}\right)= G\left(\underline{q}^n \right) + \mathcal{O}(\Delta t^{p+1}), 
\end{equation}
where $p$ is the order of accuracy of the employed RK method \cite{calvoPreservationInvariantsExplicit2006}. 

RK schemes can be made nonlinear invariant preserving by using projection techniques. With these, at each time step the solution approximated by the base integration method will be projected to the nonlinear invariant manifold $\mathcal{M}=  \biggl\{\underline{q} \; | \; G\left(\underline{q} \right)= G\left(\underline{q}_0\right)\biggr\}$. Therefore we solve the modified problem \cite{hairerGeometricNumericalIntegration2006} 
\begin{subequations} \label{projection_eq}
\begin{align}
 &\underline{\hat{q}}^{n+1} =\underline{q}^{n+1} + \lambda_n \underline{\Phi}_n, \\
&  G\left(\underline{\hat{q}}^{n+1}\right) = G\left(\underline{q}^{n}\right) , \label{Nonlin_proj_eq}
\end{align}
\end{subequations}
where $\underline{\Phi}_n$ is a search vector and $\lambda_n$ is the projection parameter to be calculated by solving the nonlinear Eq. (\ref{Nonlin_proj_eq}) \cite{hairerGeometricNumericalIntegration2006, calvoPreservationInvariantsExplicit2006}. This will be nonlinear invariant preserving, with the primary additional cost being solving a single nonlinear equation at each time step.

However, projection methods can alter certain properties of the base RK method, such as order of accuracy and the preservation of linear invariants, depending on the choice of the projection search direction \cite{calvoPreservationInvariantsExplicit2006, calvoRungeKuttaProjection2015, kojimaInvariantsPreservingSchemes2016}. While projection methods are applicable to explicit, implicit \cite{ranochaGeneralRelaxationMethods2020}, and Implicit-Explicit (IMEX) RK schemes \cite{biswasMultipleRelaxationRungeKutta2023}, this work focuses on explicit methods. In the following subsection we review existing projection methods for explicit RK methods. 

\subsection{Review of existing projection techniques}

With the standard orthogonal projection technique \cite[section IV.4]{hairerGeometricNumericalIntegration2006} the goal is to find a nonlinear invariant preserving solution that is closest to the base RK prediction. If the distance is measured in $L_2$ norm, the search vector becomes $\nabla G(\underline{\hat{q}}^{n+1})$, but due to  implicitness it is commonly replaced by the approximation
\begin{equation}
 \underline{\Phi}_{n,o}= \nabla G\left(\underline{q}^{n+1} \right) .
\end{equation}
This search direction preserves the order of accuracy of the original RK method, because it results in a sufficiently small projection correction. However, it may break preservation of linear invariants \cite{calvoPreservationInvariantsExplicit2006}, making this projection method unsuitable for the solution of hyperbolic conservation laws \cite{ranochaRelaxationRungeKutta2020}.  

Another projection direction was introduced by Del Buono and Mastroserio \cite{delbuonoExplicitMethodsBased2002} for fourth-order RK schemes with four stages to preserve inner-products $G\left(\underline{q} \right)= \underline{q}^T\underline{\underline{S}}\underline{q}= \biggl\langle \underline{q}, \underline{q} \biggr\rangle$, with $\underline{\underline{S}}$ being a constant symmetric matrix. This search direction is created by connecting the current solution to the next step solution by the original RK method, $ \underline{q}^{n+1} -\underline{q}^n$, resulting in the search vector $\underline{\Phi}_{n,r}$ 
\begin{equation}\label{relaxation_vector}
\underline{\Phi}_{n,r}= \Delta t \sum_{j=1}^s b_j\underline{R}_{j}  . 
\end{equation}
With this projection direction, linear invariants are preserved as the search direction is a linear combination of stage derivatives. However, this direction, even for inner-product invariants, causes automatic order of accuracy reduction \cite{calvoPreservationInvariantsExplicit2006, ketchesonRelaxationRungeKutta2019} . That is because it can be shown the normalized search vector becomes orthogonal to the local invariant gradient $\nabla G\left(\underline{q}^{n+1} \right)$ in the limit when time step size goes to zero
\begin{equation}
\lim_{\Delta t \to 0} \nabla G\left(\underline{q}^{n+1}\right)^T \frac{\underline{\Phi}_{n,r}}{\bigg\|\underline{\Phi}_{n,r}\bigg\|_2} =0,
\end{equation}
where $ || \cdot ||_2$ is the $L_2$ norm. A projected solution with this direction falls far from the RK prediction $\underline{q}^{n+1}$ ruining accuracy of original RK method. In \cite{delbuonoExplicitMethodsBased2002}, Del Buono and and Mastroserio showed that order of accuracy can be preserved if the projected solution $\underline{\hat{q}}^{n+1}$ is precieved as an approximation for $\underline{q}\left(t_n + \gamma_n \Delta t \right)$ instead of $\underline{q} \left(t_n + \Delta t \right)$, $\gamma_n$ being the projection parameter for the search vector \ref{relaxation_vector}. In this way the effective step size will depend on the size of projection correction.

This idea was further developed under the name of a \textit{relaxation} technique to show that it can be applied to any explicit RK method with an order of accuracy of at least two to preserve inner-products, or in general convex invariants, without order of accuracy reduction \cite{ketchesonRelaxationRungeKutta2019, ranochaRelaxationRungeKutta2020}. Moreover, For general non-convex invariants, relaxation methods retain order of accuracy if the search direction doesn't become orthogonal to the local invariant gradient faster than $\mathcal{O}(\Delta t)$ \cite{ranochaGeneralRelaxationMethods2020}. However this method may reduce efficiency when the step size is not sufficiently small. This is because $\gamma_n$ may become small, and the effective step size becomes smaller than the input step size $\gamma_n \Delta t < \Delta t$.

Calvo et al. \cite{calvoPreservationInvariantsExplicit2006} introduced an alternative projection technique known as the \textit{directional} projection method. This approach suggests using the next-step solutions from the base RK method of order at least two, alongside a lower order embedded RK method of order at least one to define the unit search vector $\underline{\Phi}_{n,d}$
\begin{equation}
 \underline{\Phi}_{n,d} = \frac{\underline{q}^{n+1} - \tilde{\underline{q}}^{n+1}}{\bigg\|\underline{q}^{n+1} - \tilde{\underline{q}}^{n+1}\bigg\|_2}  = \frac{\sum_{j=1}^s \left(b_j - \tilde{b}_j\right)\underline{R}_j}{\bigg\|\sum_{j=1}^s \left(b_j - \tilde{b}_j \right)\underline{R}_j\bigg\|_2}  ,
\end{equation}
where $\tilde{\underline{q}}^{n+1}$ represents the next-step solution by a lower order embedded method with the weights $\tilde{b}_i$. With this technique linear invariants will be preserved and the order of accuracy will depend on the choice of embedded method. Order of accuracy will be maintained if the resulting search direction doesn't become asymptotically orthogonal to $\nabla G\left(\underline{q}^{n+1}\right)$. On the other hand, if the search direction becomes nearly orthogonal to $\nabla G\left(\underline{q}^{n+1}\right)$, it leads to an order of accuracy reduction and may also cause instability by finding a projected solution that is distant from the baseline RK solution \cite{kojimaInvariantsPreservingSchemes2016}. The challenge with this method is that, except for some specific periodic problems, it may not be straightforward to predict which embedded method is suitable for a given problem \cite{calvoRungeKuttaProjection2015}.

\begin{figure}
\centering
\includegraphics[width=0.6\textwidth]{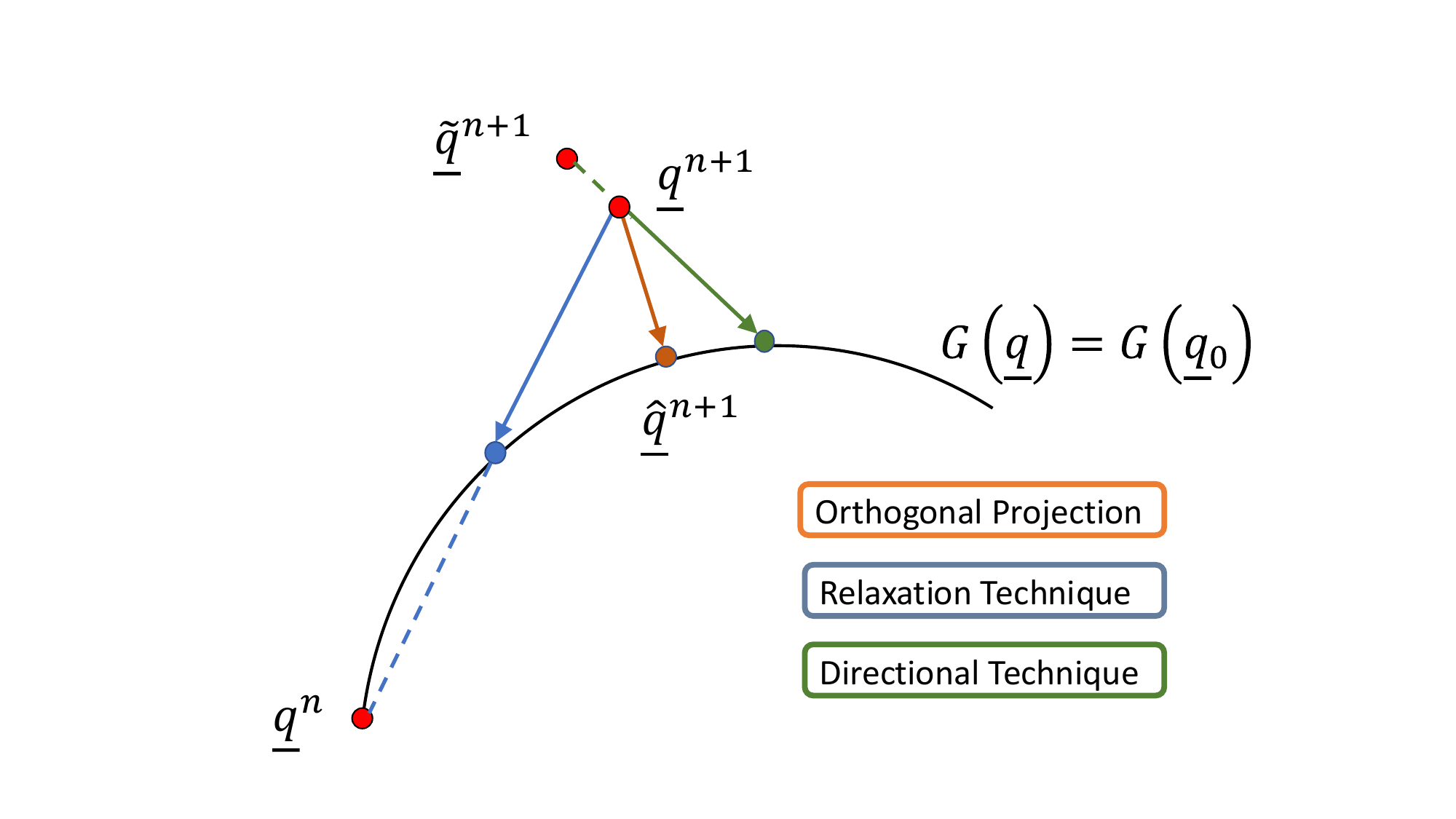}

\caption{Visualization of how each of the projection techniques project RK solution to the nonlinear invariant manifold}
\label{Figure_Projection_visualization}
\end{figure}

The direction of projection for the methods discussed above are demonstrated in Figure \ref{Figure_Projection_visualization}. From the perspective of accuracy and stability, orthogonal projection is an optimal choice, provided that preservation of linear invariants is not a concern. However, when linear invariants need to be preserved, the directional projection technique can provide search directions that preserve linear invariants and are closer to $\nabla G\left(\underline{q}^{n+1} \right)$ compared to the relaxation direction. This allows retaining accuracy without the need for step size relaxation. However, the directional projection technique requires selection of embedded RK methods, and it may not be clear beforehand which embedded method gives the best search direction for general problems, particularly at large time steps.

Projection methods can be extended to preserve multiple nonlinear invariants by employing multiple linearly independent search directions. Both the directional and the relaxation methods utilize lower-order embedded RK methods to create more than one search directions \cite{calvoPreservationInvariantsExplicit2006, biswasMultipleRelaxationRungeKutta2023}. With these methods, accuracy and stability of the projected RK will depend on the selection of embedded methods. The challenge is that it will not be clear for the user which set of embedded methods can preserve order of accuracy and is suitable for the problem under study.

Therefore, existing projection methods that preserve linear and nonlinear invariants without possibly decreasing accuracy require either relaxation of time step size, selection of lower-order embedded RK methods, or both. There is a lack of knowledge in how to obtain automatically search direction(s) that can preserve order of accuracy without need to step size relaxation.

\section{Quasi-orthogonal projection method}

In this section we present a novel and systematic strategy to define a search direction for a projection method with base RK methods of order 2 or higher. 

First, lets define the subspace spanned by a set of vectors $ \bigg\{ \underline{v}_1, ..., \underline{v}_n  \bigg\} \in \mathbb{R}^m$ denoted by $\mathit{span}\left( \underline{v}_1, ..., \underline{v}_n \right)$
\begin{equation}
\mathit{span}\left( \underline{v}_1, ..., \underline{v}_n \right)= \bigg\{ \alpha_1 \underline{v}_1+ ...+ \alpha_n \underline{v}_1 | \alpha_1, ..., \alpha_n \in \mathbb{R} \bigg\}.
\end{equation}
At each time step, we call the subspace spanned by stage derivative vectors by $\mathcal{S}$ 
\begin{equation}
\mathcal{S}:=  \mathit{span}\left(\underline{R}_1, ..., \underline{R}_s \right).
\end{equation}
We can find a set of orthonormal basis vectors $ \bigg\{ \underline{n}_1, ..., \underline{n}_k \bigg\} \in \mathbb{R}^m, \; k \leq s$, from the set of stage derivative vectors $\bigg\{ \underline{R}_1 , ..., \underline{R}_s  \bigg\}$ such that these orthonormal vectors span the subspace $\mathcal{S}$

\begin{subequations}\label{eq_ni}
\begin{align}
& \mathit{span}\left( \underline{n}_1 ..., \underline{n}_k \right)= \mathcal{S} , \\
& \underline{n}_i^T \underline{n}_j =0 , \; \mathrm{if } \; i \neq j  , \\
& \underline{n}_i^T \underline{n}_i = 1  .
\end{align}
\end{subequations}
Therefore, any linear combination of basis vectors $\underline{n}_i$ can be expressed by a linear combination stage derivative vectors $\underline{R}_i$, and vice versa.

\begin{remark}
Orthonormal basis vectors can be obtained by one by one removing linearly dependent components from stage derivative vectors and normalizing the resulting non-zero vectors.
\end{remark}

Orthogonal decomposition can be performed on $\nabla G\left( \underline{q}^{n+1}\right)$ to uniquely separate it into two components: $\nabla G_s$, which lies within subspace $\mathcal{S}$, and $\nabla G_n$, which is normal to each vector $\underline{n}_i$
\begin{subequations}
\begin{align}
& \nabla G\left(\underline{q}^{n+1}\right) = \nabla G_s + \nabla G_n  ,\\
& \nabla G_s = \sum_{i=1}^k \left( \nabla G\left(\underline{q}^{n+1}\right)^T \underline{n}_i \right) \underline{n}_i,\\
& \nabla G_s \in \mathcal{S} ,\\
&  \nabla G_n^T \underline{n}_i =0 , i= 1, .., k.
\end{align}
\end{subequations}
We propose, when $||\nabla G_s||_2$ is non-zero, using the unit vector $\nabla G_s / ||\nabla G_s||_2$ to define a method we call \textit{quasi-orthogonal} projection method
\begin{subequations}\label{Quasi_Orthogoanl_proj}
\begin{align}
& \underline{\hat{q}}^{n+1} =\underline{q}^{n+1} + \lambda_n \nabla G_s / ||\nabla G_s||_2, \\
& G\left(\underline{\hat{q}}^{n+1}\right)= G\left(\underline{q}^n\right). \label{quasi_orth_sec_eq}
\end{align}
\end{subequations}

The proposed search direction has two important properties. First, since $\nabla G_s / ||\nabla G_s||_2$ lies within the subspace $\mathcal{S}$, it can be re-created by a linear combination of stage derivatives. Therefore, this projection method preserves linear invariants and affine invariance, similar to standard RK methods. Secondly, in the following we show that projection of $\nabla G\left( \underline{q}^{n+1} \right)$ onto the proposed direction has the largest magnitude compared to projection of $\nabla G\left(\underline{q}^{n+1}\right)$ onto any unit vector $\underline{d} \in \mathcal{S}$. This property will be used in showing the order of accuracy of the proposed projection method.

\begin{lemma}\label{nabla_G_biggest_innerprod}
If there is a unit vector $\underline{d} \in \mathcal{S}$ such that $\nabla G\left(\underline{q}^{n+1} \right)^T \underline{d} \neq 0$, then $||\nabla G_s||_2 \neq 0$, and
\begin{equation}
\bigg|  \nabla G\left(\underline{q}^{n+1} \right)^T \frac{\nabla G_s}{|| \nabla G_s ||_2} \bigg|= ||\nabla G_s ||_2  \geq  \bigg|\nabla G\left(\underline{q}^{n+1} \right)^T \underline{d} \bigg|. 
\end{equation}

\end{lemma}

\begin{proof}[Proof of Lemma \ref{nabla_G_biggest_innerprod}]
We have
\begin{equation}
 \nabla G \left(\underline{q}^{n+1}\right)^T \underline{d}= \left(\nabla G_s + \nabla G_n \right)^T \underline{d} = \nabla G_s^T\underline{d} . 
\end{equation}
So, if $\nabla G\left(\underline{q}^{n+1} \right)^T \underline{d}$ is non-zero, $||\nabla G_s||_2$ should be also non-zero. Moreover, the maximum value for $| \nabla G(\underline{q}^{n+1})^T \underline{d} |$ becomes $||\nabla G_s||_2$.
\end{proof}

\begin{lemma}\label{d_for_convex_G}
Assume problem (\ref{IVP}) integrated by an RK method of order at least two is non-stationary and contains a convex invariant $G\left(\underline{q}\right)$, i.e. having $ R\left(t_n, \underline{q}^n \right)^T \nabla^2 G\left(\underline{q}^n \right)R\left(t_n, \underline{q}^n \right) >0$ at each time step. Therefore, a directional projection method using a first order embedded method results in a unit search vector $\underline{\Phi}_{n,d} \in \mathcal{S}$ such that 
\begin{equation}
\lim_{\Delta t \rightarrow 0} \nabla G\left(\underline{q}^{n+1}\right)^T \underline{\Phi}_{n,d} \neq 0 
\end{equation}

\end{lemma}

\begin{proof}[Proof of Lemma \ref{d_for_convex_G} ]
For the directional projection method \cite{calvoPreservationInvariantsExplicit2006} we know that using a base RK method of order $p \geq 2$ and an embedded RK method of order $q$ ($ 1 \leq q \leq p-1$) gives a search direction $ \underline{\Phi}_{n,d} =  \underline{C}_{q+1}/||\underline{C}_{q+1}||_2 + \mathcal{O}(\Delta t)$, where $\underline{C}_{q+1}\Delta t^{q+1}$ is the leading error term of the embedded method. Therefore, using a first order embedded method gives the unit search direction
\begin{equation}
 \underline{\Phi}_{n,d}= \frac{\underline{C}_2}{\Big\|\underline{C}_2\Big\|_2} + \mathcal{O}(\Delta t),
\end{equation}
where we have \cite{hairerSolvingOrdinaryDifferential1996}
\begin{equation}\label{C_2}
\frac{\underline{C}_2}{ \Big\| \underline{C}_2 \Big\|_2}= \frac{  \left(   R_t + R_{\underline{q}} \; R \right)(t_n, \underline{q}^n)  }{\bigg\| \left(   R_t + R_{\underline{q}} \; R \right)(t_n, \underline{q}^n)    \bigg\|_2 } .
\end{equation}
On the other hand, since $\nabla G\left(\underline{q}\right)^T R\left(\underline{q},t\right)$ is zero according to Eq. (\ref{dG_dt}) for $t \in [0, T]$ and $\underline{q} \in R^m$, its time derivative also should be zero for a conservative ODE system
\begin{equation}\label{G_second_derivative}
\frac{d}{dt}\left(\nabla G\left(\underline{q}\right)^T R\left(\underline{q},t\right)\right)= R\left(\underline{q},t\right)^T \nabla^2 G\left(\underline{q}\right) R\left(\underline{q},t\right) + \nabla G\left(\underline{q}\right)^T \left( R_t + R_{\underline{q}} \; R \right) (t, \underline{q}) =0. 
\end{equation}
At $(t_n, \underline{q}^n)$ the first term on the right-hand side is non-zero for non-stationary problems with convex invariants. Consequently, the second term must also be nonzero at $(t^n, \underline{q}^n)$
\begin{equation}
\nabla G\left(\underline{q}^n\right)^T \left( R_t + R_{\underline{q}} R \right) (t^n, \underline{q}^n) \neq 0. 
\end{equation}
As $\nabla G\left( \underline{q}^{n+1}\right)= \nabla G\left(\underline{q}^n\right) + \mathcal{O}(\Delta t)$, we have	
\begin{equation}
 \lim_{\Delta t \rightarrow 0} \nabla G\left(\underline{q}^{n+1}\right)^T \underline{\Phi}_{n,d} \neq 0 ,
\end{equation}
and since $ \underline{\Phi}_{n,d} \in \mathcal{S}$, proof is complete.
\end{proof}

From the proof of Lemma \ref{d_for_convex_G} one direct result is the following Corollary.

\begin{corollary}\label{accuracy_directional_convex}
A directional projection method \cite{calvoPreservationInvariantsExplicit2006} using first order embedded methods for non-stationary problems containing a  convex invariant is solvable and preserves order of accuracy for sufficiently small step sizes.

\end{corollary}

\begin{proof}[Proof of Corollary \ref{accuracy_directional_convex}]

According to \cite[Theorem 4.1.]{calvoPreservationInvariantsExplicit2006}, a directional projection method with a first order embedded method is solvable and preserves accuracy for small enough step sizes if $\nabla G\left(\underline{q}^n \right)^T \underline{C}_2 \neq 0$, which is shown to be true for non-stationary problems with a convex invariant in the proof of Lemma \ref{d_for_convex_G}.

\end{proof}

\begin{remark}
As mentioned also in \cite{calvoRungeKuttaProjection2015}, all first order embedded methods in a directional projection method yield asymptotically identical search direction as $\Delta t \rightarrow 0$. However, when the step size is relatively large, it is not always clear in advance which first order embedded method will be most suitable for the problem at hand \cite{kojimaInvariantsPreservingSchemes2016}. 

\end{remark}

Now we can show solvability and accuracy of the proposed projection method for convex invariants.
\begin{theorem}\label{Theorem_convexG}
For non-stationary problems (\ref{IVP}) with a convex invariant integrated with RK methods of order $p\geq 2$, there exists $h^* >0$ such that for $\Delta t\in (0, h^*]$ we have $||\nabla G_s||_2 \neq 0$ and the projection method (\ref{Quasi_Orthogoanl_proj}) is solvable and the resulting order of accuracy will be $\hat{p} \geq p$.
\end{theorem}

\begin{proof}[Proof of Theorem \ref{Theorem_convexG}]
According to Lemma \ref{d_for_convex_G}, with the prior assumptions there exists a unit vector $\underline{d} \in \mathcal{S}$ such that for sufficiently small step sizes $ \lim_{\Delta t \rightarrow 0} \nabla G \left(\underline{q}^{n+1} \right)^T \underline{d} \neq 0$. Therefore, according to Lemma \ref{nabla_G_biggest_innerprod}, $||\nabla G_s||_2 $ will be nonzero for sufficiently small step sizes.

The reminder of proof closely follows the approach outlined in \cite[Theorem 4.1.]{calvoPreservationInvariantsExplicit2006}. We define a real smooth function $g(\lambda , \Delta t), \; \lambda \in  \mathbb{R}$
\begin{equation}
g(\lambda , \Delta t)= G \left( \underline{q}^{n+1} + \lambda \frac{\nabla G_s}{||\nabla G_s||_2} \right) - G(\underline{q}^n).
\end{equation}
We have $g(0,0)=0$, and 
\begin{equation}
\frac{\partial g}{\partial \lambda} (0,0)= \lim_{\Delta t \rightarrow 0} \nabla G(\underline{q}^{n+1})^T \frac{\nabla G_s}{||\nabla G_s||_2}= \lim_{\Delta t \rightarrow 0} ||\nabla G_s||_2 \neq 0 .
\end{equation}
Therefore, solvability of the proposed projection method for sufficiently small step sizes will be proven using the implicit function theorem, and the resulting integration method will have an order of accuracy of $\hat{p} \geq p$.
\end{proof}

\begin{remark}
For inner-product invariants after obtaining $\Phi_n = \nabla G_s/ ||\nabla G_s||_2$ at each time step, Eq. (\ref{quasi_orth_sec_eq}) becomes a quadratic equation and can be solved analytically similar to \cite[Eq. 18]{calvoPreservationInvariantsExplicit2006}. Moreover, it can be shown that, for sufficiently small step sizes, $\nabla G_s/ ||\nabla G_s||_2$ results in the smallest projection correction among all search directions in $\mathcal{S}$.
\end{remark}

Then, we continue this section by showing order of accuracy and solvability of the proposed projection method for general invariants compared to relaxation and directional projection methods.

\begin{theorem}\label{Theorem_generalG_relaxation}
Assume that for an ODE system (\ref{IVP}) with a general (not-necessarily convex) invariant function, a relaxation method for sufficiently small step sizes is solvable and preserves accuracy by ensuring $\nabla G\left(\underline{q}^{n+1}\right)^T (\underline{q}^{n+1} - \underline{q}^n)/||\underline{q}^{n+1}-\underline{q}^n||_2= c\Delta t + \mathcal{O}(\Delta t^2), \; c\neq 0$ on non-constant step sizes \cite{ranochaGeneralRelaxationMethods2020}. Then for small enough step sizes $||\nabla G_s||_2$ is non-zero and the proposed projection method (\ref{Quasi_Orthogoanl_proj}) is solvable and preserves accuracy with a constant step size. 
\end{theorem}

\begin{proof}[Proof of Theorem \ref{Theorem_generalG_relaxation}]
Here parameter $c$ can be obtained and it is
\begin{equation}
 c=
\frac{1}{||\sum_{i=1}^sb_i \underline{R}_i||_2} \left(R\left(t_n, \underline{q}^n \right)^T \nabla^2 G\left(\underline{q}^n\right)R\left(t_n, \underline{q}^n \right) + \frac{1}{2} \nabla G\left(\underline{q}^n\right)\left( R_t + R_{\underline{q}}R \right) (t_n, \underline{q}^n) \right).
\end{equation}
Having $c\neq 0$ alongside Eq. (\ref{G_second_derivative}) at $(t_n, \underline{q}^n)$ results in
\begin{equation}
 \nabla G\left(\underline{q}^n \right)^T\left( R_t + R_{\underline{q}}R \right) (t_n, \underline{q}^n) \neq 0.
\end{equation}
Therefore, for this problem the directional projection method with a first order embedded RK method gives a unit search direction $ \underline{\Phi}_{n,d} =  \underline{C}_{q+1}/||\underline{C}_{q+1}||_2 + \mathcal{O}(\Delta t)$ such that $\nabla G\left(\underline{q}^{n+1}\right)^T \underline{\Phi}_{n,d}$ becomes non-zero for sufficiently small step sizes. and the rest of proof is similar to the proof for Theorem \ref{Theorem_convexG}.
\end{proof}

\begin{theorem}\label{Theorem_generalG_directional}
Assume that for the ODE systems (\ref{IVP}) with a general invariant function a directional projection method with a base RK method of order $p\geq 2$ results in a unit search direction $\underline{\Phi}_{n,d}$ such that $\nabla G\left( \underline{q}^{n+1}\right)^T \underline{\Phi}_{n,d} = B(\underline{q}^n) \Delta t^r + \mathcal{O}(\Delta t^{r+1})$ with $0 \leq r \leq p/2$ and $B(\underline{q}^n) \neq 0$, guaranteeing an order of accuracy of at least least $p-r$. Then for sufficiently small step sizes $||\nabla G_s||_2$ is non-zero and the proposed projection method (\ref{Quasi_Orthogoanl_proj}) is solvable and guarantees an order of accuracy of $\hat{p} \geq p-r$.
\end{theorem}

\begin{proof}[Proof of Theorem \ref{Theorem_generalG_directional}]
Since $\underline{\Phi}_{n,d}\in \mathcal{S}$ and with sufficiently small step sizes $\nabla G\left( \underline{q}^{n+1}\right)^T \underline{\Phi}_{n,d} \neq 0$, according to Lemma \ref{nabla_G_biggest_innerprod} for small step sizes $||\nabla G_s||_2 \neq 0$ and
\begin{equation}
\bigg|\nabla G\left(\underline{q}^{n+1}\right)^T \frac{\nabla G_s}{||\nabla G_s||_2}\bigg|= ||\nabla G_s||_2 \geq \bigg|\nabla G\left( \underline{q}^{n+1}\right)^T \underline{\Phi}_{n,d}\bigg| .
\end{equation}
The rest of proof is similar to corresponding part of \cite[Theorem 4.1.]{calvoPreservationInvariantsExplicit2006}.
\end{proof}

Therefore, the quasi-orthogonal projection method provides a systematic strategy to determine the projection search direction. The provable order of accuracy with this method will be at least as high as the order of accuracy with the relaxation and directional projection methods. However, the proposed method doesn't require a change in the step size or selection of embedded RK methods, alleviating the limitations of relaxation and directional projection methods.

Regarding computational cost, with the proposed method typically the main cost would be solving a single nonlinear equation (using iterative methods for general invariant functions). Therefore, the total computational cost with the quasi-orthogonal projection method remains minimal compared to base RK methods.

\begin{remark}
For small systems, subspace $\mathcal{S}$ may span the entire $\mathbb{R}^m$ domain. In such cases, we would have $\nabla G_s= \nabla G\left( \underline{q}^{n+1}\right)$, making the proposed projection method equivalent to the standard orthogonal projection method.
\end{remark}

\begin{remark}

There might be situations that the search vector $\nabla G_s$ becomes zero in the limit when the step size goes to zero. For stationary problems, stage derivatives become zero and consequently the vector $\nabla G_s$ becomes zero. In such cases, the invariant function remains constant and there would be no need for projection correction. However, if for a non-stationary problem with a general invariant function $\nabla G_s$ becomes zero for sufficiently small step sizes, according to Lemma \ref{nabla_G_biggest_innerprod}, for any unit search vector $\underline{d} \in \mathcal{S}$ we would have $\nabla G\left(\underline{q}^{n+1} \right)^T \underline{d} = 0$. Therefore, any search direction in $\mathcal{S}$ can cause a severe projection correction, which in turn may lead to lose of accuracy and instability issues.

\end{remark}

\subsection{Extension to dissipative systems}

We call ODE system (\ref{IVP}) dissipative if its function $G\left(\underline{q}\right)$ is dissipative in time
\begin{equation}\label{dissipative_G}
 \frac{d}{dt}G\left(\underline{q}\left(t\right)\right)= \nabla G\left(\underline{q}\right)^TR\left(t, \underline{q}\right) < 0 ,  \qquad \textrm{for all} \; t \in [0, T],
\end{equation}
and preserving monotonicity by a time integration method requires 
\begin{equation}
G\left(\underline{q}^{n+1}\right) < G\left(\underline{q}^n\right). 
\end{equation}
However, many explicit RK schemes don't guarantee monotonicity preservation for general problems, even with small step sizes \cite{ketchesonRelaxationRungeKutta2019, ranochaRelaxationRungeKutta2020, sunStabilityFourthOrder2017}, which can potentially lead to nonphysical solutions. The projection technique can mitigate this shortcoming by making RK schemes monotonicity preserving using a proper projection equation, as demonstrated in \cite{ketchesonRelaxationRungeKutta2019, ranochaRelaxationRungeKutta2020}.

Now we extend quasi-orthogonal projection method to dissipative systems. We propose, when $||\nabla G_s||_2$ is non-zero, using the search direction $ \nabla G_s / ||\nabla G_s||_2$ with the following modified problem 
\begin{subequations} \label{projection_eqs_dissipative}
\begin{align}
& \underline{\hat{q}}^{n+1} =\underline{q}^{n+1} + \lambda_n \nabla G_s / ||\nabla G_s||_2, \\
&  G\left(\underline{\hat{q}}^{n+1}\right) = G\left(\underline{q}^{n}\right) + \Delta t \sum_{i=1}^s  \nabla G\left(\underline{q}_i\right)^T \underline{R}_i, \label{Nonlin_proj_eq_dissipative}
\end{align}
\end{subequations}
where Eq. (\ref{Nonlin_proj_eq_dissipative}) guarantees monotonicity preservation for the projected solution. The following theorem demonstrates solvability and accuracy of the proposed projection method for dissipative systems.

\begin{theorem}\label{theorem_qo_dissipative}
For dissipative ODE systems integrated with explicit RK methods of order $p\geq 2$, there exists $h^* >0$ such that for $\Delta t \in (0, h^*]$ we have $||\nabla G_s||_2 \neq 0$ and projection method (\ref{projection_eqs_dissipative}) is solvable and the resulting order of accuracy will be $\hat{p} \geq p$.
\end{theorem}

\begin{proof}[Proof of Theorem \ref{theorem_qo_dissipative}]
We have
\begin{equation}
\nabla G\left(\underline{q}^{n+1}\right)^T R\left(t_n, \underline{q}^n \right) = \nabla G\left(\underline{q}^{n}\right)^T R\left(t_n, \underline{q}^n \right) + \mathcal{O}(\Delta t)  . 
\end{equation}
Since $\nabla G\left(\underline{q}^{n}\right)^T R\left(t_n, \underline{q}^n \right) <0$, for a small step size $\nabla G\left(\underline{q}^{n+1}\right)^T R\left(t_n, \underline{q}^n \right) /||R\left(t_n, \underline{q}^n\right)||_2 \neq 0$. So, according to
Lemma \ref{nabla_G_biggest_innerprod} for small enough step sizes $||\nabla G_s||_2 \neq 0$. 

Then, we define the smooth function $g(\Delta t, \lambda), \; \lambda \in \mathbb{R}$
\begin{equation}
g(\lambda , \Delta t)= G\left( \underline{q}^{n+1} + \lambda \frac{\nabla G_s}{||\nabla G_s||_2}\right) - G\left(\underline{q}^n\right) - \Delta t  \sum_{i=1}^s  \nabla G\left(\underline{q}_i\right)^T \underline{R}_i,
\end{equation}
and we have
\begin{equation}
g(0,0)=0,
\end{equation}
\begin{equation}
\frac{ \partial g}{\partial\lambda}(0,0)= \lim_{\Delta t \rightarrow 0} \nabla G(\underline{q}^{n+1})^T\frac{\nabla G_s}{||\nabla G_s||_2}= ||\nabla G_s||_2 \neq 0 .
\end{equation}
So, as before similar to \cite[Theorem 4.1.]{calvoPreservationInvariantsExplicit2006} solvability of Eq. \ref{projection_eqs_dissipative} for small enough step size can be demonstrated using the implicit function theorem. 

Moreover, according to \cite[Corollary 2.13.]{ranochaRelaxationRungeKutta2020} we have $g(0, \Delta t)=\mathcal{O}(\Delta t^{p+1})$. Additionally we have 
\begin{equation}
\frac{\partial g}{\partial \lambda}(0, \Delta t)= \frac{\partial g}{\partial \lambda}(0, 0) + \mathcal{O}(\Delta t). 
\end{equation}
Therefore, a projection parameter $\lambda_n$ satisfying $g(\lambda_n, \Delta t)=0$ will be $\lambda_n = \mathcal{O}(\Delta^{p+1})$ and the resulting projection method will have an order of accuracy of $\hat{p} \geq p$ \cite{calvoPreservationInvariantsExplicit2006}.

\end{proof}

\subsection{Extension to multiple invariants}

Suppose now ODE system (\ref{IVP}) has $\ell \geq 1$ smooth invariant functions $G_1\left(\underline{q}\right), ..., G_{\ell}\left(\underline{q}\right)$
\begin{equation}
\frac{d}{dt} G_j\left(\underline{q}(t)\right)=  \nabla G_j\left(\underline{q}\right)^T R\left(t,\underline{q}\right)= 0, \quad j=1,..., \ell,  \quad \textrm{for all} \; t \in [0, T] .
\end{equation}
A projection technique with explicit RK methods can be used to preserve the function $G= \left( G_1, ..., G_{\ell} \right)^T: \mathbb{R}^m \rightarrow \mathbb{R}^{\ell}$ at each time step, by defining $\ell$ linearly independent search directions.

With a directional projection method, Calvo et al. \cite{calvoPreservationInvariantsExplicit2006} proposed using $\ell$ linearly independent lower order RK methods to define the search directions. The relaxation method also has been extended to preserve multiple invariants in \cite{biswasMultipleRelaxationRungeKutta2023} by relaxing the step size and using $\ell -1$ lower order embedded methods. However, with both methods the user has to select a set of embedded RK methods, and this selection defines the order of accuracy of the projection method. 

For multiple invariants we propose performing orthogonal decomposition at each time step on $\nabla G_j\left(\underline{q}^{n+1}\right)$, $j= 1, ..., \ell$, using orthonormal vectors $ \bigg\{ \underline{n}_1, ..., \underline{n}_k \bigg\} $ defined in Eq. (\ref{eq_ni}) to uniquely find $\nabla G_{j,s}$

\begin{subequations} 
\begin{align}
 & \nabla G_j\left(\underline{q}^{n+1}\right) = \nabla G_{j,s} + \nabla G_{j,n} , \\
& \nabla G_{j,s} = \sum_{i=1}^k \left(  \nabla G_j\left(\underline{q}^{n+1}\right)^T \underline{n}_i \right) \underline{n}_i, \\
&  \nabla G_{j,s} \in \mathcal{S} , \\
&  \nabla G_{j,n}^T \; \underline{n}_i =0 , i= 1, .., k.
\end{align}
\end{subequations}
When $||\nabla G_{j,s}||_2 \neq 0$ for $j=1, ..., \ell$, we suggest using search directions $\nabla G_{j,s} / ||\nabla G_{j,s}||_2$ to define a quasi-orthogonal projection method for multiple invariants
\begin{subequations}\label{Multiple_qo}
\begin{align}
& \underline{\hat{q}}^{n+1} =\underline{q}^{n+1} + \sum_{j=1}^{\ell} \lambda_n^{(j)} \nabla G_{j,s} / ||\nabla G_{j,s}||_2, \\
& G\left(\underline{\hat{q}}^{n+1}\right)= G\left(\underline{q}^n\right), \label{Multiple_qo_eq}
\end{align}
\end{subequations}
where $\lambda_n^{(j)} \in  \mathbb{R}$ are projection parameters to be obtained to satisfy Eq. (\ref{Multiple_qo_eq}).

To proceed, for a set of unit vectors $\bigg\{\underline{d}_1, ..., \underline{d}_{\ell}\bigg\} \in \mathcal{S}$ we define the matrix $M\left(\underline{d}_1, ..., \underline{d}_{\ell}\right)$
\begin{equation}
M\left(\underline{d}_1, ..., \underline{d}_{\ell}\right) := \lim_{\Delta t \rightarrow 0}
\begin{bmatrix}
\nabla G_1\left(\underline{q}^{n+1}\right)^T\underline{d}_1 & \dots & \nabla G_1\left(\underline{q}^{n+1}\right)^T\underline{d}_{\ell} \\
\vdots & \ddots & \vdots \\
\nabla G_{\ell}\left(\underline{q}^{n+1}\right)^T\underline{d}_1 & \dots & \nabla G_{\ell}\left(\underline{q}^{n+1}\right)^T\underline{d}_{\ell}
\end{bmatrix},
\end{equation}
and for nonzero $\nabla G_{j,s}$ vectors, we define matrix $\underline{\underline{M}}_s$ as $M(\nabla G_{1,s}/||\nabla G_{1,s}||_2, ...,  \nabla G_{\ell,s}/||\nabla G_{\ell,s}||_2)$ 
\begin{equation}
\underline{\underline{M}}_s := \lim_{\Delta t \rightarrow 0}
\begin{bmatrix}
\nabla G_1\left(\underline{q}^{n+1}\right)^T\nabla G_{1,s}/||\nabla G_{1,s}||_2  & \dots & \nabla G_1\left(\underline{q}^{n+1}\right)^T\nabla G_{\ell,s}/||\nabla G_{\ell,s}||_2 \\
\vdots & \ddots & \vdots \\
\nabla G_{\ell}\left(\underline{q}^{n+1}\right)^T\nabla G_{1,s}/||\nabla G_{1,s}||_2 & \dots & \nabla G_{\ell}\left(\underline{q}^{n+1}\right)^T\nabla G_{\ell,s}/||\nabla G_{\ell,s}||_2
\end{bmatrix}.
\end{equation}

\begin{lemma}\label{lemma_nonsingular_matrices}
If there is a set of unit vectors $\bigg\{\underline{d}_1, ..., \underline{d}_{\ell}\bigg\} \in \mathcal{S}$ such that matrix $M\left(\underline{d}_1, ..., \underline{d}_{\ell}\right)$ becomes nonsingular, then for sufficiently small step sizes we have $||\nabla G_{j,s}||_2\neq 0$ for $j=1, .., \ell$, and matrix $\underline{\underline{M}}_s$ is also nonsingular.
\end{lemma}

\begin{proof}[Proof of Lemma \ref{lemma_nonsingular_matrices}]
Since the $\underline{d}_j$ vectors are in $\mathcal{S}$, the matrix $M\left(\underline{d}_1, ..., \underline{d}_{\ell}\right)$ can be written as
\begin{equation}
M\left(\underline{d}_1, ..., \underline{d}_{\ell}\right)= \lim_{\Delta t \rightarrow 0}
\begin{bmatrix}
\nabla G_{1,s}^T\underline{d}_1 & \dots & \nabla G_{1,s}^T\underline{d}_{\ell} \\
\vdots & \ddots & \vdots \\
\nabla G_{\ell,s}^T\underline{d}_1 & \dots & \nabla G_{\ell,s}^T\underline{d}_{\ell}
\end{bmatrix}.
\end{equation}
If this matrix is non-singular, its rows (and columns) should be linearly independent. Therefore, for small enough step sizes vectors $\nabla G_{j,s}$ should be linearly independent and nonzero: $||\nabla G_{j,s}|| \neq 0$. 

Moreover, vectors 
\begin{equation}
\begin{bmatrix}
\nabla G_{1,s}^T\nabla G_{j,s}  \\
\vdots  \\
\nabla G_{\ell,s}^T\nabla G_{j,s}
\end{bmatrix}, \; j=1, ..., \ell,
\end{equation}
can be shown to be linearly independent: given linear independence of vectors $\nabla G_{j,s}$, the only way to get 
\begin{equation}
\begin{bmatrix}
\nabla G_{1,s}^T \left(
\sum_{j=1}^{\ell} k_j \nabla G_{j,s} \right)  \\
\vdots  \\
\nabla G_{\ell,s}^T \left(
\sum_{j=1}^{\ell} k_j \nabla G_{j,s} \right)
\end{bmatrix} = \underline{0}
\end{equation}
with scalar multipliers $k_j$ is to have the trivial solution, $k_j =0$ for $j=1, ..., \ell$. Therefore, the columns and rows of matrix $\underline{\underline{M}}_s$ are linearly independent, and matrix $\underline{\underline{M}}_s$ is nonsingular.
\end{proof}

We can now demonstrate solvability and order of accuracy of proposed projection method for multiple invariants.

\begin{theorem}\label{theorem_multiple_qo}
Assume that ODE system (\ref{IVP}) has $\ell$ invariant functions $G_1, ..., G_{\ell}$, and a directional projection method with a base RK method of order $p\geq 2$ can preserve invariants using unit search directions $\underline{\Phi}_{n,d}^{(1)}, ..., \underline{\Phi}_{n,d}^{(\ell)}$ such that the matrix $M\left(\underline{\Phi}_{n,d}^{(1)}, ..., \underline{\Phi}_{n,d}^{(\ell)}\right)$ is nonsingular and consequently the order of accuracy is preserved \cite{calvoPreservationInvariantsExplicit2006}. Therefore, for sufficiently small step sizes vectors $\nabla G_{s,1}, .., \nabla G_{s,\ell}$ will be nonzero and the quasi-orthogonal method for multiple invariants (\ref{Multiple_qo}) is solvable and the resulting order of accuracy will be $\hat{p} \geq p$.
\end{theorem}

\begin{proof}[Proof of Theorem \ref{theorem_multiple_qo}]
According to Lemma \ref{lemma_nonsingular_matrices} since matrix $M\left(\underline{\Phi}_{n,d}^{(1)}, ..., \underline{\Phi}_{n,d}^{(\ell)}\right)$ is non-singular, for small enough step sizes vectors $\nabla G_{s,1}, .., \nabla G_{s,\ell}$ are nonzero and matrix $\underline{\underline{M}}_s$ is nonsingular. Moreover, similar to the proof provided in the related part of \cite[Theorem 4.2.]{calvoPreservationInvariantsExplicit2006} we can define $\underline{\lambda}_n= \left( \lambda_n^{(1)}, ...,   \lambda_n^{(\ell)} \right)^T$ and function $g(\underline{\lambda}_n, \Delta t): \mathbb{R}^{\ell} \times \mathbb{R} \rightarrow \mathbb{R}^{\ell}$
\begin{equation}
g(\underline{\lambda}_n, \Delta t)= G \Biggl( \underline{q}^{n+1} + \sum_{j=1}^{\ell} \lambda_n^{(j)} \nabla G_{j,s} / ||\nabla G_{j,s}||_2 \Biggr) - G \left(  \underline{q}^{n} \right) . 
\end{equation}
In a way similar to the proof provided for Theorem \ref{Theorem_convexG} we can show that for sufficiently small step sizes Eq. (\ref{Multiple_qo}) is solvable and the order of accuracy after projection will be $\hat{q} \geq p$.
\end{proof}

\begin{remark}
Theorem \ref{theorem_multiple_qo} demonstrates that the proposed method systematically employs an optimal set of search directions to preserve the order of accuracy. This is advantageous compared to both directional projection and relaxation methods, which require selection of embedded RK methods by the user in the process of creating search directions. This selection may automatically result in a decrease in the order of accuracy and stability issues.
\end{remark}

\begin{remark}
It can be shown that to preserve accuracy with a quasi-orthogonal approach (and the directional projection method) for $\ell$ invariants, a necessary condition is to have a base RK method with number of stages $s \geq \ell +1$.
\end{remark}







\section{Numerical examples} \label{Sec_exs}

This section presents several numerical experiments to examine the proposed quasi-orthogonal projection method and to compare it with existing projection techniques. To best perform a thorough comparison, selected example cases from papers \cite{ketchesonRelaxationRungeKutta2019, ranochaRelaxationRungeKutta2020, biswasMultipleRelaxationRungeKutta2023} are reproduced here. These examples involve a range of linear and non-linear illustrative problems using the following base explicit RK schemes:
\begin{itemize}
\item SSPRK(2,2): Second-order, two-stage SSP method of \cite{shu_efficient_1988}. 
\item RK(3,3): Third-order, three-stage standard RK scheme, named RK31 in \cite{butcherNumericalMethodsOrdinary2008}.
\item Heun(3,3): Third-order, three-stage method of Heun, can be found in appendix of \cite{biswasMultipleRelaxationRungeKutta2023}.
\item RK(4,4): Classical fourth-order RK scheme with four stages, named RK41 in \cite{butcherNumericalMethodsOrdinary2008}. 
\item DP(7,5): Fifth-order, seven-stage method of \cite{princeHighOrderEmbedded1981}.
\item BSRK(8,5): Fifth-order, eight-stage method of \cite{bogackiEfficientRungeKuttaPair1996}.
\end{itemize}

\subsection{Linear dissipative system}

Sun \& Shu \cite{sunStabilityFourthOrder2017} demonstrated that for a semidiscrete system which is dissipative with respect to the squared norm, or energy, time integration with RK(4,4) can adversely cause an increase in energy after the first integration step, regardless of how small the step size is. As provided in \cite{ketchesonRelaxationRungeKutta2019}, an indicative example is the linear dissipative system in the form of $\frac{d}{dt}\underline{q}\left(t\right)=\underline{\underline{L}} \; \underline{q}\left(t\right)$

\begin{equation}\label{Eq_lin_decay}
\frac{d}{dt}
\begin{bmatrix}
q_1  \\
q_2 \\
q_3
\end{bmatrix}= \begin{bmatrix}
-1 & -2 & -2 \\
0 & -1 & -2 \\
0 & 0 & -1
\end{bmatrix} \begin{bmatrix}
q_1 \\
q_2 \\
q_3
\end{bmatrix},
\end{equation}
with an initial condition equal to the first right singular vector of $\mathcal{R}(0.5\underline{\underline{L}} )$, with $\mathcal{R}(z)$ being the stability polynomial of RK(4,4). For this problem time rate change of energy is negative
\[ \frac{d}{dt} \bigg\|\underline{q}(t)\bigg\|_2^2  \leq 0.  \]

\begin{figure}[h]
\centering
\includegraphics[width=0.4\textwidth]{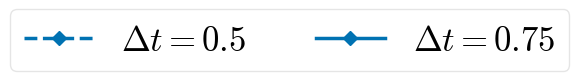}

\begin{subfigure}[b]{0.3\textwidth}
\centering
\includegraphics[width=\textwidth]{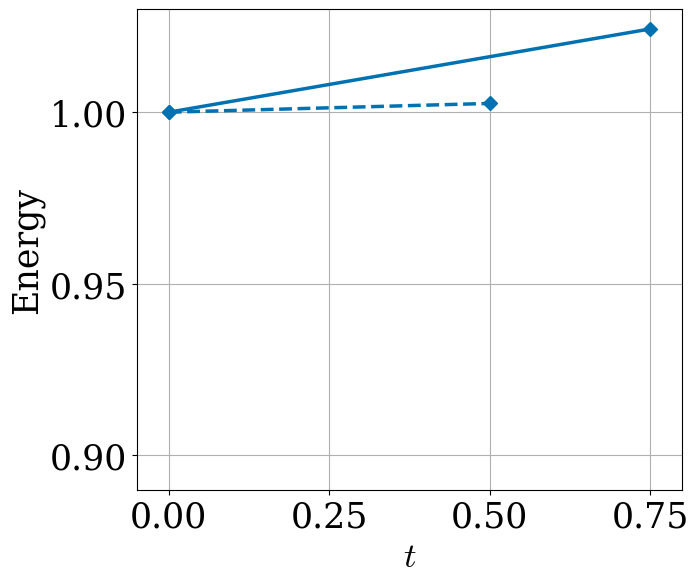}
\caption{ RK(4,4) }
\label{Figure_3_evol_a}
\end{subfigure}
\begin{subfigure}[b]{0.3\textwidth}
\centering
\includegraphics[width=\textwidth]{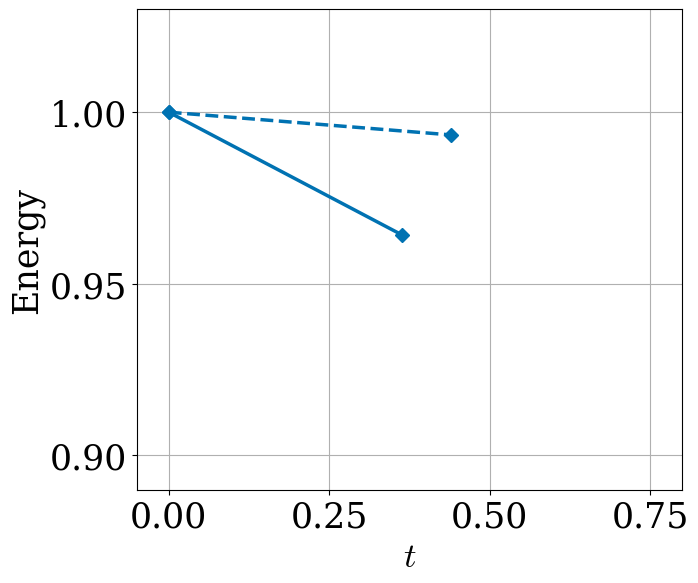}
\caption{Relaxation method}
\label{Figure_3_evol_b}
\end{subfigure}
\begin{subfigure}[b]{0.3\textwidth}
\centering
\includegraphics[width=\textwidth]{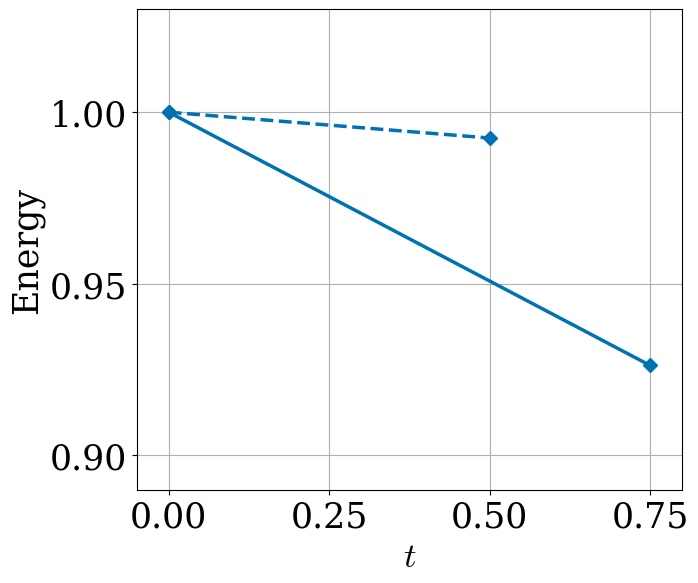}
\caption{Quasi-orthogonal method}
\label{Figure_3_evol_c}
\end{subfigure}

\caption{ Energy change after the first time step for energy dissipative ODE system (\ref{Eq_lin_decay}) integrated with standard RK(4,4), relaxation method with RK(4,4), and quasi-orthogonal projection method with RK(4,4). }
\label{Fig_linear_energy_decay_compare}

\end{figure}

Figure \ref{Fig_linear_energy_decay_compare} shows the change in energy after the first time step for ODE system (\ref{Eq_lin_decay}) integrated with standard RK(4,4), relaxation RK(4,4), and quasi-orthogonal projection with RK(4,4). At both tested step sizes, $\Delta t= 0.5$ and $\Delta t= 0.7$, RK(4,4) causes energy to increase after the first time step. However, both relaxation and quas-orthogonal projection techniques preserve monotonicity of problem by ensuring a decrease in energy. This figure also shows that for both time step sizes, the actual time step size with the relaxation method is reduced to below $0.5$.

\begin{figure}[h]
\centering

\includegraphics[width=0.5\textwidth]{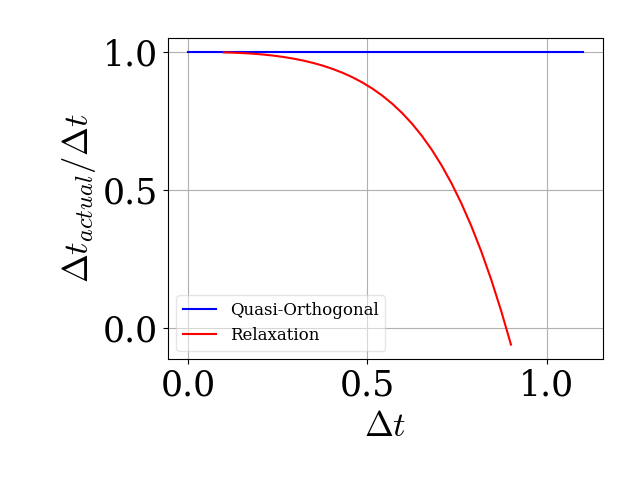}
\caption{ Actual time step size over input step size for relaxation and quasi-orthogonal projection methods after the first time step integrating ODE system (\ref{Eq_lin_decay}).}

\label{Fig_linear_energy_decay_step_size}

\end{figure}

Moreover, Figure \ref{Fig_linear_energy_decay_step_size} shows the change in the actual step size indicated by $\Delta t_{actual}$ over input step size after the first time step for relaxation and quasi-orthogonal projection methods while integrating the ODE system (\ref{Eq_lin_decay}) with a base method of RK(4,4). With the relaxation method $\Delta t_{actual}/\Delta t$ monotonically decreases while increasing the input step size, until it reaches zero at around $\Delta t \approx 0.89$. On the other hand, the proposed quasi-orthogonal method keeps the time step size unchanged and it provides a larger solvability region (up to $\Delta t \approx 1.1$).

\subsection{Nonlinear oscillator}

Following another example from the work of Ketcheson \cite{ketchesonRelaxationRungeKutta2019}, we test integration schemes on the nonlinear oscillator problem
\begin{equation}\label{Eq_nonlin_oscil}
\frac{d}{dt}
\begin{bmatrix}
q_1 \\
q_2
\end{bmatrix} = \frac{1}{\bigg\|\underline{q}\bigg\|_2^2}\begin{bmatrix}
-q_2 \\
q_1
\end{bmatrix}, \; \begin{bmatrix}
q_1(0) \\
q_2(0)
\end{bmatrix}= \begin{bmatrix}
1\\
0
\end{bmatrix}.
\end{equation}
The time rate change of energy for this problem is zero 
\[ \frac{d}{dt} \bigg\| \underline{q}(t) \bigg\|_2^2 =0, \]
and it has the following exact analytical solution
\[ \begin{bmatrix}
q_1(t) \\
q_2(t)
\end{bmatrix} = \begin{bmatrix}
\text{cos}(t)\\
\text{sin}(t) \\
\end{bmatrix}. \]

\begin{figure}
\centering
\includegraphics[width=\textwidth]{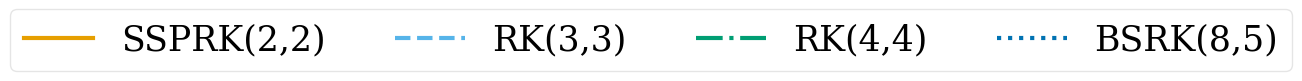}

\begin{subfigure}[b]{0.3\textwidth}
\centering
\includegraphics[width=\textwidth]{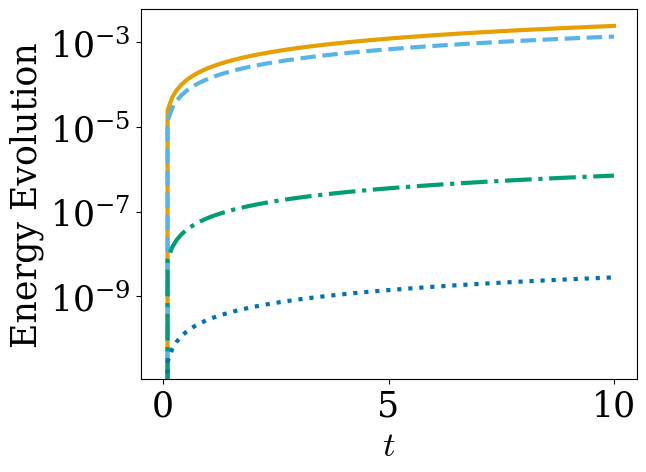}
\caption{ Base RK methods }
\label{Figure_3_evol_a}
\end{subfigure}
\begin{subfigure}[b]{0.3\textwidth}
\centering
\includegraphics[width=\textwidth]{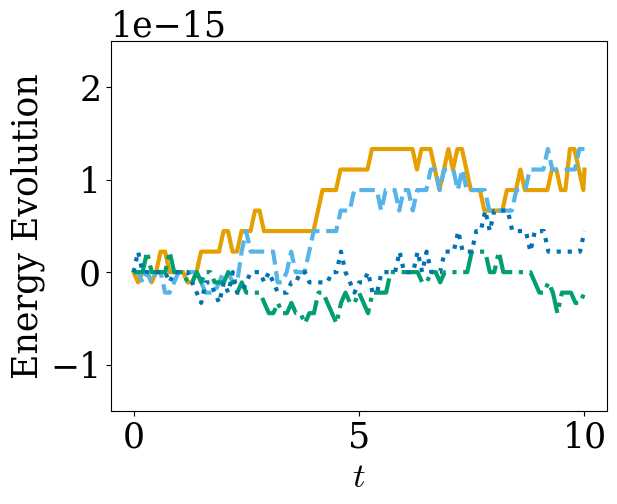}
\caption{Relaxation methods}
\label{Figure_3_evol_b}
\end{subfigure}
\begin{subfigure}[b]{0.3\textwidth}
\centering
\includegraphics[width=\textwidth]{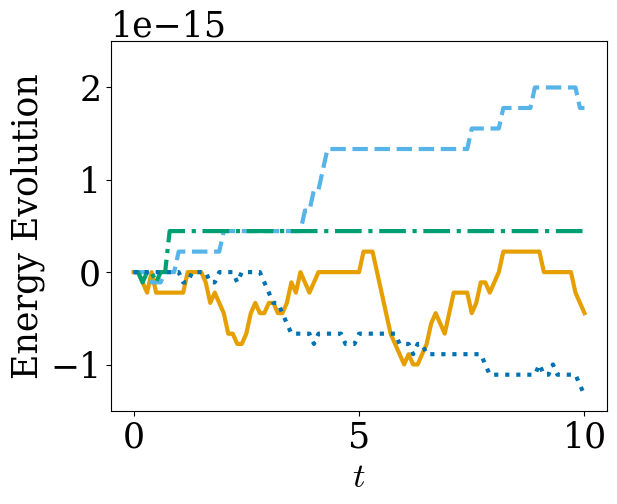}
\caption{Quasi-orthogonal methods}
\label{Figure_3_evol_c}
\end{subfigure}

\caption{Energy evolution in the nonlinear oscillator problem  (\ref{Eq_nonlin_oscil}) integrated with different integration schemes and a time step of $\Delta t=0.1$. All the original RK methods change energy up to their truncation error, but with the relaxation and quasi-orthogonal projection methods energy is preserved up to machine precision. }

\label{Figure_Nonlin_oscil_evol}

\end{figure}

Figure \ref{Figure_Nonlin_oscil_evol} shows that employing each of the unmodified RK schemes with a time step size of $\Delta t=0.1$ leads to a monotonic increase in energy. On the other hand, both relaxation and quasi-orthogonal projection methods ensure the base RK schemes conserve energy up to machine precision.

\begin{figure}
\centering
\includegraphics[width=0.6\textwidth]{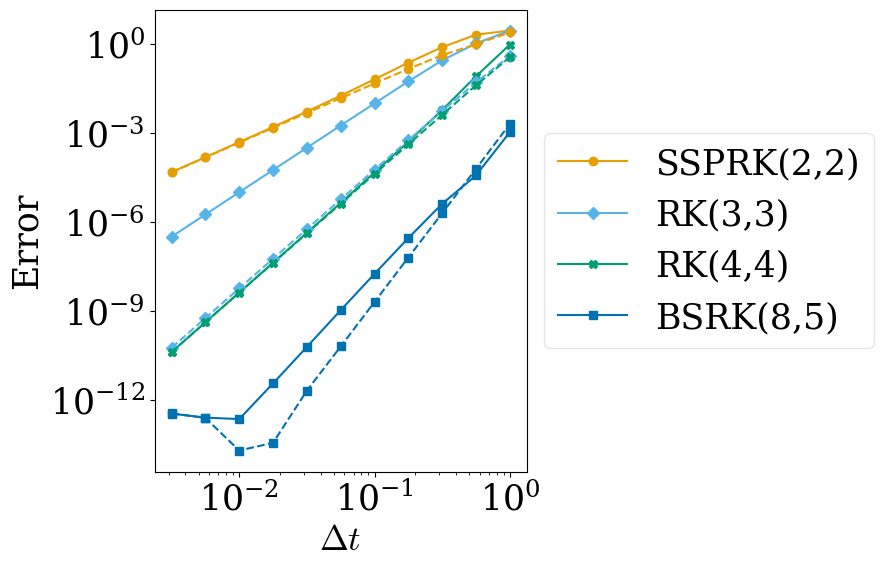}

\caption{Convergence study for the nonlinear oscillator problem (\ref{Eq_nonlin_oscil}). For base RK methods convergence is demonstrated by solid lines, while for their quasi-orthogonal projection counterpart it is depicted by dashed lines. With the proposed method, order of accuracy is equal to, or higher than, the corresponding base RK scheme. }
\label{Figure_3_convergence}
\end{figure}

Regarding the convergence rates, Figure \ref{Figure_3_convergence} shows solution convergence at time $t=10$ for unmodified schemes with solid lines and the corresponding invariant-preserving schemes with the quasi-orthogonal projection method in dashed lines. It confirms that with the proposed method order of accuracy is equal to, or higher than, the base RK method for a non-stationary problem containing a convex invariant.

\begin{figure}
\centering
\includegraphics[width=0.8\textwidth]{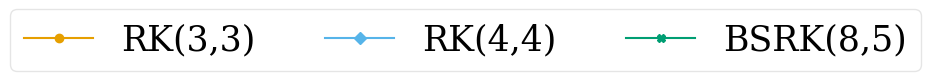}

\begin{subfigure}[b]{0.3\textwidth}
\centering
\includegraphics[width=\textwidth]{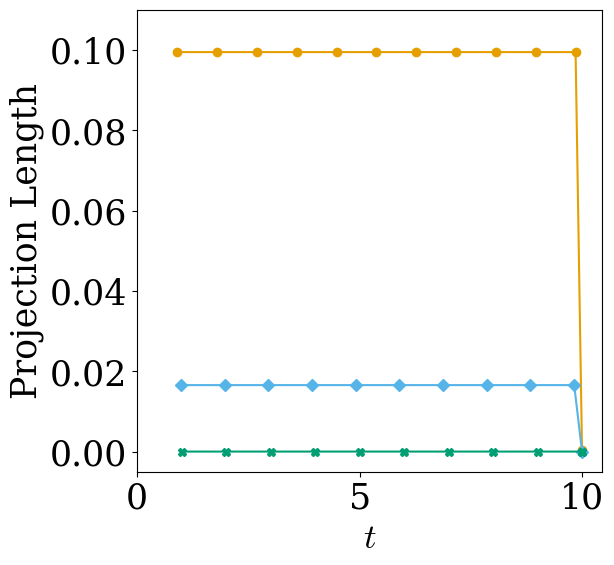}
\caption{ Relaxation methods }
\label{Figure_3_evol_a}
\end{subfigure}
\begin{subfigure}[b]{0.3\textwidth}
\centering
\includegraphics[width=\textwidth]{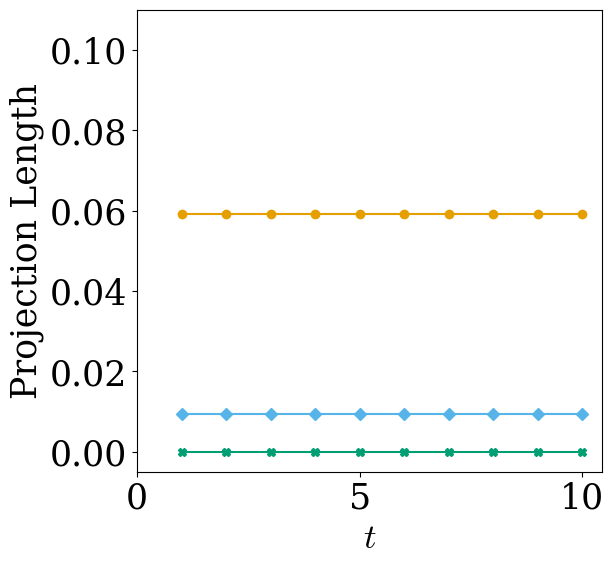}
\caption{Directional projection methods}
\label{Figure_3_evol_b}
\end{subfigure}
\begin{subfigure}[b]{0.3\textwidth}
\centering
\includegraphics[width=\textwidth]{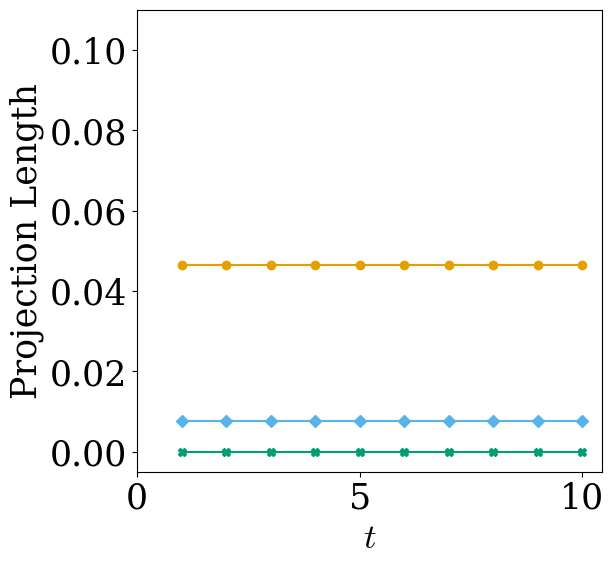}
\caption{Quasi-orthogonal methods}
\label{Figure_3_evol_c}
\end{subfigure}

\caption{Projection length using different projection schemes for the nonlinear oscillator problem (\ref{Eq_nonlin_oscil}). }
\label{Figure_3_projection_length}
\end{figure}

For this problem with different projection methods we can compare the projection length
\[ \bigg\|\underline{\hat{q}}^{n+1} -\underline{q}^{n+1} \bigg\|_2 . \]
Figure \ref{Figure_3_projection_length} compares projection length with $\Delta t=1$ for the relaxation method, the directional projection method with the first order Euler method as the embedded method, and the proposed quasi-orthogonal projection method. It shows that the quasi-orthogonal projection method among other mentioned methods leads to minimum projection length for an inner-product norm invariant.

With this problem, since stage derivatives cover the whole 2D domain, the quasi-orthogonal projection method will behave identical to the orthogonal projection technique.








\subsection{Burgers equation}

Another example case from \cite{ketchesonRelaxationRungeKutta2019} is the inviscid Burger's problem with the PDE 
\begin{equation}
\frac{d}{dt}Q+ \frac{1}{2}\frac{d}{dx}(Q^2)=0  ,
\end{equation}
on a periodic interval of $-1 \leq x \leq 1$ with the initial condition
\[ Q(x,0)= e^{-30x^2}.  \]

\subsubsection{Burgers energy conservative}

This problem can be transformed to an energy conservative ODE system by discretizing the domain with $50$ equally-spaced points and using the second-order accurate symmetric flux \cite{tadmorEntropyStabilityTheory2003}
\begin{equation}
q_i'(t)= -\frac{1}{\Delta x} (F_{i+1/2} - F_{i-1/2}) , \qquad F_{i+1/2}=\frac{q_i^2 + q_iq_{i+1} + q_{i+1}^2}{6}.
\end{equation}

\begin{figure}[t] 
\centering
\includegraphics[width=0.7\textwidth]{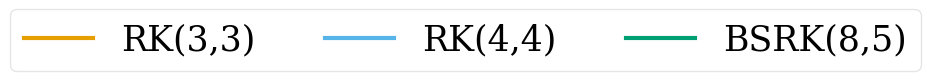}

\begin{subfigure}[b]{0.3\textwidth}
\centering
\includegraphics[width=\textwidth]{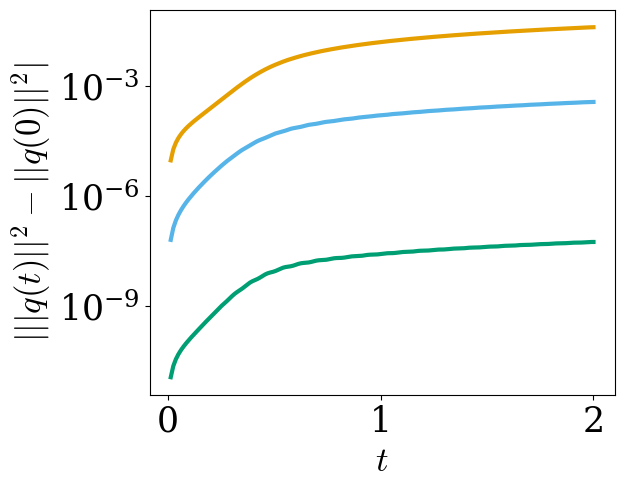}
\caption{ Base RK methods }
\end{subfigure}
\begin{subfigure}[b]{0.3\textwidth}
\centering
\includegraphics[width=\textwidth]{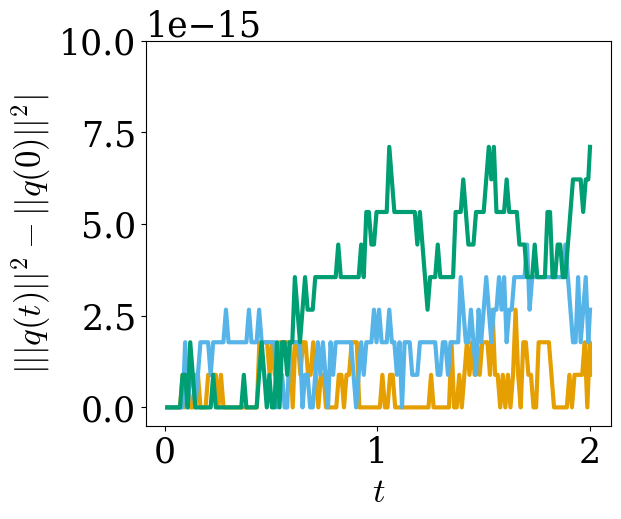}
\caption{Relaxation methods}
\end{subfigure}
\begin{subfigure}[b]{0.3\textwidth}
\centering
\includegraphics[width=\textwidth]{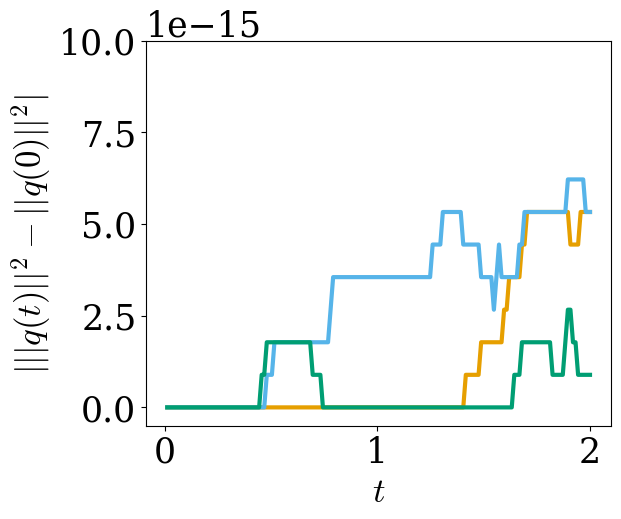}
\caption{Quasi-orthogonal methods}
\end{subfigure}

\caption{Evolution of energy for the Burgers' equation with different integration schemes with a step size of $\Delta t=0.3\Delta x$ up to a final time of $t_f=2$. With each of the unmodified RK schemes energy increases monotonically, while invariant preserving counterparts conserve energy up to machine precision}
\label{Figure_Burgers_evol}
\end{figure}

Figure \ref{Figure_Burgers_evol} shows energy change in the system while integrating with a time step size of $\Delta t=0.3\Delta x$ with unmodified RK methods and invariant conservative counterparts. It demonstrates that while with base RK methods energy changes with time, both the relaxation and the quasi-orthogonal projection methods keep energy unchanged over time.

\begin{figure}[t] 
\centering
\includegraphics[width=0.7\textwidth]{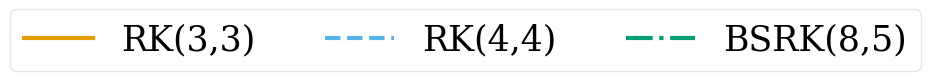}

\begin{subfigure}[b]{0.3\textwidth}
\centering
\includegraphics[width=\textwidth]{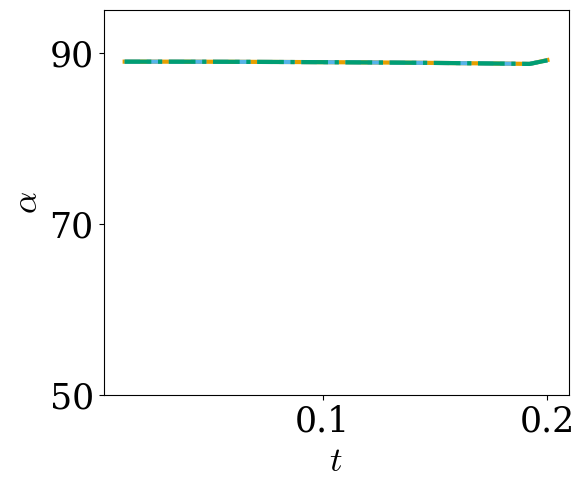}
\caption{ Relaxation methods }
\end{subfigure}
\begin{subfigure}[b]{0.3\textwidth}
\centering
\includegraphics[width=\textwidth]{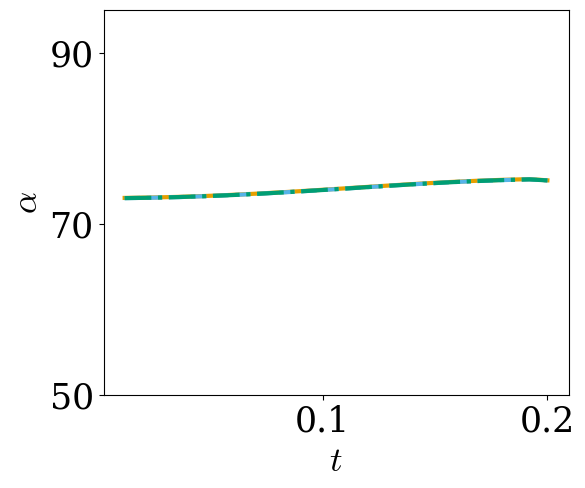}
\caption{Directional projection methods}
\end{subfigure}
\begin{subfigure}[b]{0.3\textwidth}
\centering
\includegraphics[width=\textwidth]{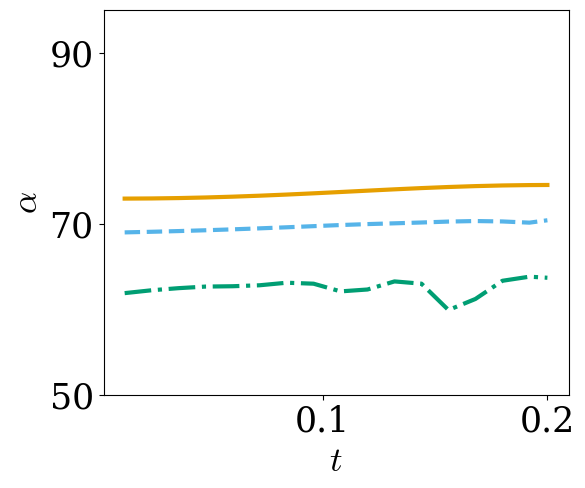}
\caption{Quasi-orthogonal methods}
\end{subfigure}

\caption{Comparison of angles ($\alpha$) in degrees created by search directions and $\nabla G\left(\underline{q}^{n+1}\right)$ by each of the projection methods.}
\label{Figure_Burgers_proj_angle}
\end{figure}

We can also compare search directions employed by the relaxation technique, a directional projection method with the first order Euler method as the embedded method, and the proposed quasi-orthogonal method in terms of angle (in degrees) created by $\nabla G\left( \underline{q}^{n+1}\right)$ and the search direction $\underline{\Phi}_n$
\[ \alpha=  \cos^{-1} \left( \frac{\nabla G\left(\underline{q}^{n+1}\right)^T \underline{\Phi}_n}{||\nabla G\left(\underline{q}^{n+1}\right)||_2 \;||\underline{\Phi}_n||_2}  \right) .\]
Figure \ref{Figure_Burgers_proj_angle} demonstrates evolution of $\alpha$ over time for each of the invariant preserving methods. It shows that while the relaxation method results in a search direction almost orthogonal to $\nabla G \left(\underline{q}^{n+1}\right)$, the directional projection method and the proposed method lead to search directions closer to $\nabla G \left(\underline{q}^{n+1}\right)$. Moreover, by increasing the number of stages of the base RK method, the proposed projection technique achieved smaller values for $\alpha$.

\begin{figure}[t] 
\centering
\includegraphics[width=0.7\textwidth]{Burgers_evol_legend.png}

\includegraphics[width=0.4\textwidth]{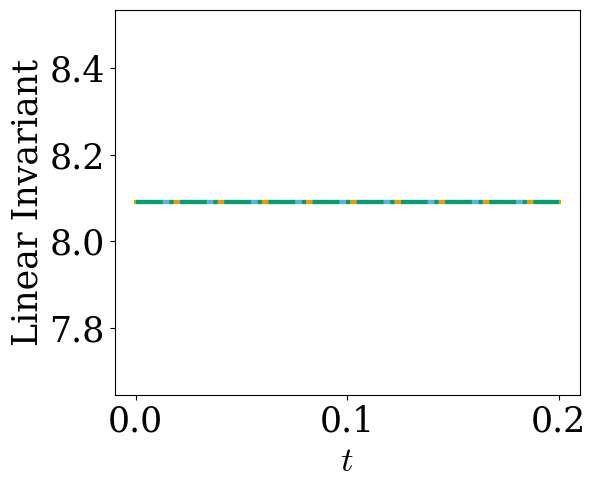}

\caption{Change in linear invariant of the Burgers' problem with the quasi-orthogonal projection method. }
\label{Figure_Burgers_linear_invariant}
\end{figure}

In Figure \ref{Figure_Burgers_linear_invariant} we demonstrate that with each of base RK methods, the proposed projection method preserved the linear invariant of the problem, which is $\underline{1}^T \underline{q}$, $\underline{1}$ being vector of ones.

\begin{figure}[t] 
\centering

\includegraphics[width=0.7\textwidth]{Burgers_evol_legend.png}

\includegraphics[width=0.3\textwidth]{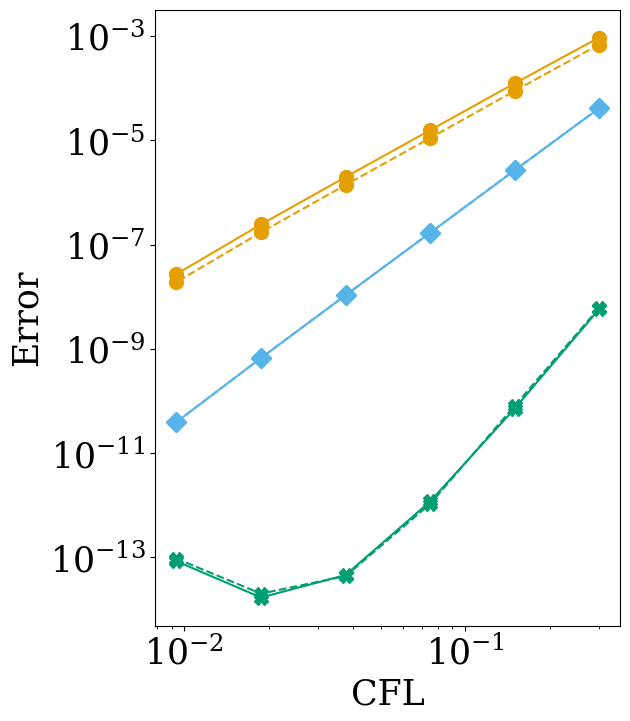}

\caption{Convergence analysis for Burger's equation integrated up to a final time of $t_f=0.2$ with the time steps: $\Delta t= \text{CFL} \times \Delta x$. Solid lines represent base integration methods and dashed lines represent quasi-orthogonal projection counterparts. It confirms that the quasi-orthogonal projection method preserves order of accuracy of the base RK methods for convex invariants.}
\label{Fig_Burgers_convergence}
\end{figure}

Then, to visualize order of accuracy of the quasi-orthogonal projection method, convergence analysis has been performed at $t_f=0.2$ with a fixed spatial discretization and different input time steps: $\Delta t= \text{CFL} \times \Delta x$, $\text{CFL}= 0.3 \times 0.5^{0, 1, ..,6}$. Results in Figure \ref{Fig_Burgers_convergence} confirm that orders of accuracy of base RK methods are preserved with the quasi-orthogonal projection method.






\subsection{Rigid Body Rotation}

With this example we examine behavior of the proposed projection method applied to Euler equations
\begin{subequations}\label{Eq_RBR}
\begin{align}
&  \frac{d}{dt}q_1 = (\alpha - \beta)q_2 q_3 ,\\
&  \frac{d}{dt}q_2 = (1-\alpha)q_3 q_1 ,\\
&  \frac{d}{dt}q_3 = (\beta -  1)q_1 q_2 ,
\end{align}
\end{subequations}
with $\left(q_1(0), q_2(0), q_3(0)\right)^T= (0,1,1)^T$, where $\alpha= 1 + 1/\sqrt{1.51}$ and $\beta= 1- 0.51/\sqrt{1.51}$. This ODE system describes motion of a free rigid body with its center of mass at the origin in terms of its angular momenta \cite{calvoPreservationInvariantsExplicit2006, biswasMultipleRelaxationRungeKutta2023}. There are two invariant functions for this problem
\begin{subequations}
\begin{align}
&  G_1(q_1, q_2, q_3) = q_1^2 + q_2^2 + q_3^2 ,\\
&  G_2(q_1, q_2, q_3)= q_1^2 + \beta q_2^2 + \alpha q_3^2 .
\end{align}
\end{subequations}

\begin{figure}[t] 
\centering
\includegraphics[width=0.8\textwidth]{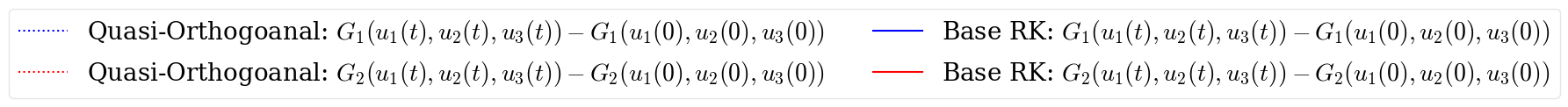}

\begin{subfigure}[b]{0.3\textwidth}
\centering
\includegraphics[width=\textwidth]{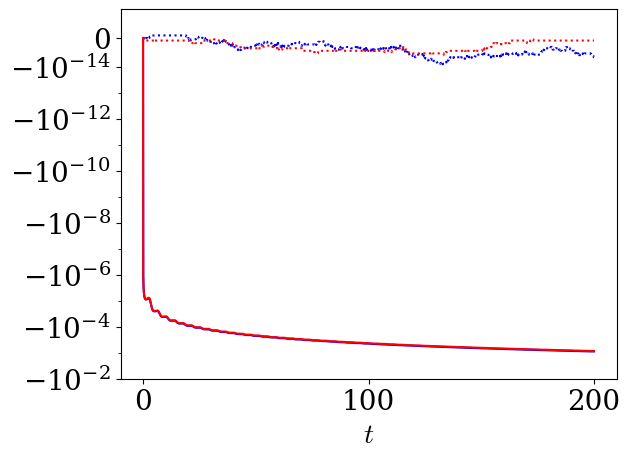}
\caption{ Heun(3,3) }
\end{subfigure}
\begin{subfigure}[b]{0.3\textwidth}
\centering
\includegraphics[width=\textwidth]{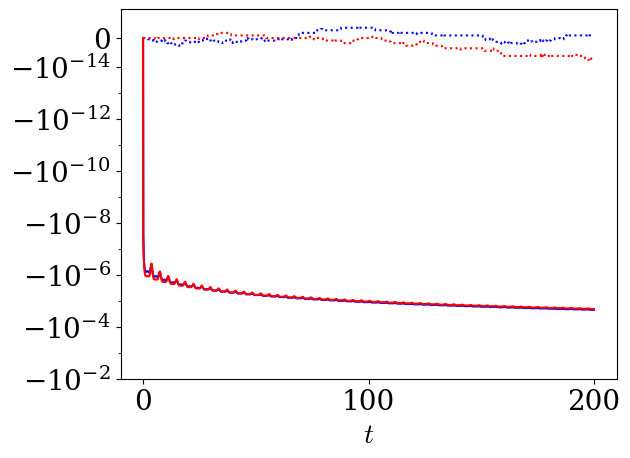}
\caption{ RK(4,4) }
\end{subfigure}
\begin{subfigure}[b]{0.3\textwidth}
\centering
\includegraphics[width=\textwidth]{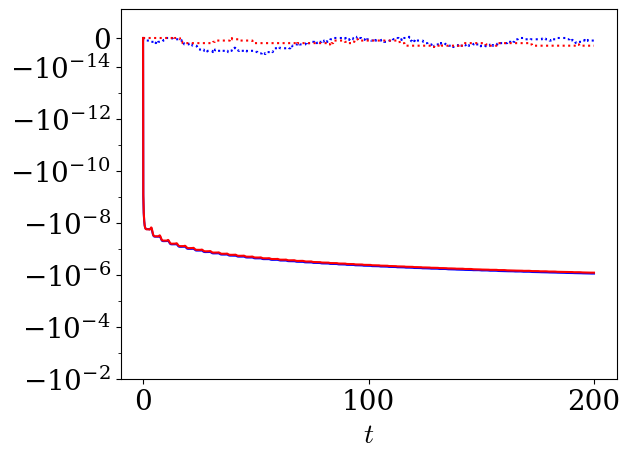}
\caption{ DP(7,5) }
\end{subfigure}

\caption{Change in invariants for rigid body rotation problem (\ref{Eq_RBR}) integrated with base RK and corresponding quasi-orthogonal projection methods. }
\label{Figure_RBR_invar_change}
\end{figure}

Similar to \cite{biswasMultipleRelaxationRungeKutta2023}, we examine change in invariants with and without employment of invariant-preserving methods for Heun(3,3) with $\Delta t=0.04$, RK(4,4) with $\Delta t=0.1$, and DP(7,5) with $\Delta t=0.1$.  Figure \ref{Figure_RBR_invar_change} shows changes in invariants for the base methods and the quasi-orthogonal projection counterparts, demonstrating that the proposed method preserved the two invariants.

\begin{figure}[t] 
\centering
\includegraphics[width=0.5\textwidth]{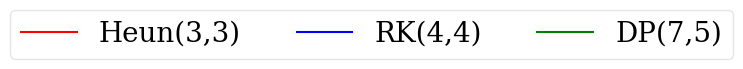}

\begin{subfigure}[b]{0.3\textwidth}
\centering
\includegraphics[width=\textwidth]{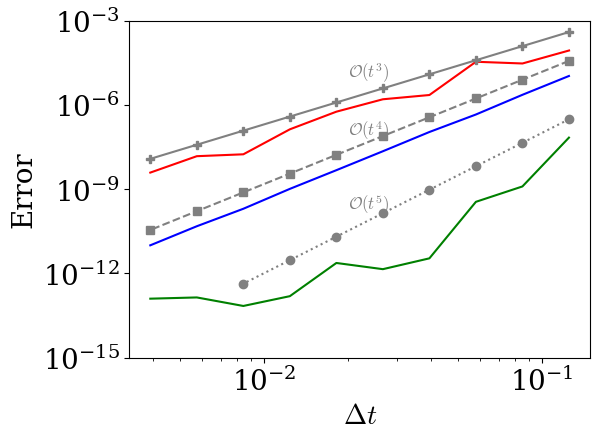}
\caption{ Relaxation methods }
\end{subfigure}
\begin{subfigure}[b]{0.3\textwidth}
\centering
\includegraphics[width=\textwidth]{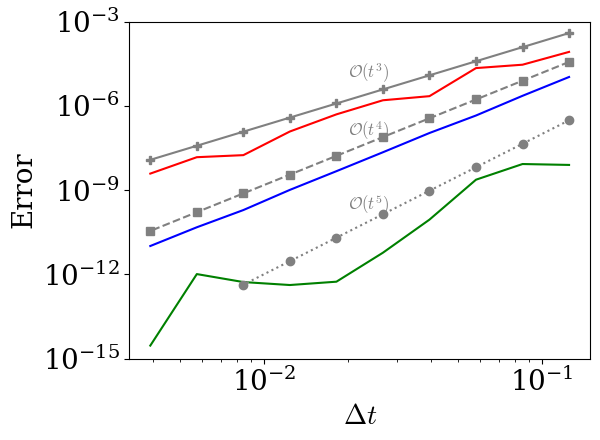}
\caption{ Directional projection methods }
\end{subfigure}
\begin{subfigure}[b]{0.3\textwidth}
\centering
\includegraphics[width=\textwidth]{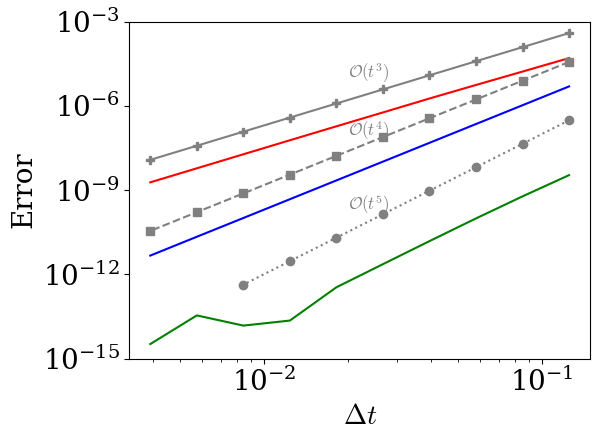}
\caption{ Quasi-Orthogonal methods }
\end{subfigure}

\caption{Convergence study for the rigid body rotation problem (\ref{Eq_RBR}) using each of the invariant-preserving RK methods.}
\label{Figure_RBR_convergence}
\end{figure}

Then, we perform convergence analysis to demonstrate convergence rates with the proposed quasi-orthogonal approach compared to relaxation and directional projection methods. For relaxation and directional projection methods, since they require embedded RK methods, the ones outlined in \cite{biswasMultipleRelaxationRungeKutta2023} along with the first order euler method (for directional projection method) are employed. Results provided in Figure \ref{Figure_RBR_convergence} show that the quasi-orthogonal method, without needing to use embedded RK methods, provides a smoother convergence rate using each of the base RK methods, compared to other invariant-preserving methods.


\section{Conclusions} \label{Sec_conclusion}

We have proposed a new family of invariant-preserving projection methods for explicit RK integration schemes. This method is applicable to ODE systems containing one or multiple invariants, as well as dissipative systems. The search direction(s) with this method are systematically obtained through an efficient approach and they provide advantageous properties over existing projection techniques. The resulting search direction(s) preserve linear invariants of the problem, similar to standard RK schemes. They are proven to be optimal in preserving the order of accuracy among all search directions created by RK stage derivative vectors, without needing to use embedded RK schemes. Moreover, this method preserves order of accuracy without the need for step size relaxation. Numerical results show that these aforementioned properties are observed in practice for a dissipative system and problems containing one or multiple invariants. Future work will focus on efficient extension of the proposed method to implicit and IMEX RK schemes.


\section*{Acknowledgments}

The authors acknowledge financial support from the Natural Sciences and Engineering Research Council of Canada (NSERC) and the Fonds de Recherche du Québec - Nature et Technologies (FRQNT) via the NOVA program. 

Moreover, the authors acknowledge that for replication of the results from papers \cite{ketchesonRelaxationRungeKutta2019, ranochaRelaxationRungeKutta2020, biswasMultipleRelaxationRungeKutta2023} and creating associated plots some codes in open access repositories provided by \cite{ranochaCodeConvexRelaxationRungeKutta2019, biswasCodeMultiplerelaxationRunge2023} are used in this work.

\bibliographystyle{unsrt}

\bibliography{Time_integration}


\end{document}